\documentclass[10pt,a4paper]{article}
\usepackage{amsmath}
\usepackage{latexsym}
\usepackage{amssymb}
\usepackage{amstext}
\usepackage{amsfonts}
\usepackage{amsthm}
\usepackage{fancyhdr}
\usepackage{indentfirst}
\usepackage{graphicx}
\usepackage{epsfig}
\usepackage{pgf}
\usepackage{mathrsfs}
\usepackage{listings}
\usepackage{cases}
\usepackage{epstopdf}
\usepackage{graphicx}
\usepackage{enumerate}
\usepackage{abstract}
\usepackage{caption}
\date{}

\begin{document}
\title{The generalized reciprocal distance matrix of graphs}
\author{Gui-Xian Tian$^a$\footnote{Corresponding author. E-mail addresses: gxtian@zjnu.cn (G-X. Tian), 2980614556@qq.com (M-J. Cheng), cuishuyu@zjnu.cn (S-Y. Cui).}, Mei-Jiao Cheng$^a$, Shu-Yu Cui$^{b}$ \\%EndAName
    {\small{\it $^a$Department of Mathematics,}}
    {\small{\it Zhejiang Normal University, Jinhua 321004, China}}\\
    {\small{\it $^b$Xingzhi College, Zhejiang Normal University, Jinhua 321004, China}}
}\maketitle
\begin{abstract}
Let $G$ be a simple undirected connected graph with the Harary matrix $RD(G)$, which is also called the reciprocal distance matrix of $G$. The reciprocal distance signless Laplacian matrix of $G$ is $RQ(G)=RT(G)+RD(G)$, where $RT(G)$ denotes the diagonal matrix of the vertex reciprocal transmissions of graph $G$. This paper intends to introduce a new matrix $RD_{\alpha}(G)=\alpha RT(G)+(1-\alpha)RD(G)$, $\alpha\in [0,1]$, to track the gradual change from $RD(G)$ to $RQ(G)$. First, we describe completely the eigenvalues of $RD_{\alpha}(G)$ of some special graphs. Then we obtain serval basic properties of $RD_{\alpha}(G)$ including inequalities that involve the spectral radii of the reciprocal distance matrix, reciprocal distance signless Laplacian matrix and $RD_{\alpha}$-matrix of $G$. We also provide some lower and upper bounds of the spectral radius of $RD_{\alpha}$-matrix. Finally, we depict the extremal graphs with maximal spectral radius of the $RD_{\alpha}$-matrix among all connected graphs of fixed order and precise vertex connectivity, edge connectivity, chromatic number and independence number, respectively.


\end{abstract}

\section{Introduction}
All the graphs studied in this paper are simple undirected connected graphs. Let $G=(V(G),E(G))$ be a connected graph of order $n$, where its vertex set is $V(G)=\{v_{1},v_{2},\ldots, v_{n}\}$ and its edge set is $E(G)$. We use $v_{i}v_{j}$ to denote that two vertices $v_{i}$ and $v_{j}$ are neighbors in the graph $G$. The adjacency matrix of $G$ is the matrix $A(G)=(a_{i,j})_{n\times n}$, where $a_{i,j}=1$ if $v_{i}v_{j}\in E(G)$ and $0$, otherwise. The degree diagonal matrix of $G$ is denote by $\overline{D}(G)$.
Then the signless Laplacian matrix and Laplacian matrix are $Q(G)=\overline{D}(G)+A(G)$ and $L(G)=\overline{D}(G)-A(G)$, respectively. The sum of distances of $v$ from all other vertices of $G$ are denote by $Tr_{G}(v)$, which is called as the vertex transmission of $v$ in $G$. The distance matrix of $G$ is denoted by $D(G)=(d_{ij})$, where $d_{ij}$ represents the distance between $v_{i}$ and $v_{j}$. Let the diagonal matrix $Tr(G)$ be the vertex transmissions matrix with the $(i,i)$th entry of $Tr(G)$ being $Tr_{G}(v_{i})$.
Similarly,  the distance signless Laplacian matrix and distance Laplacian matrix of graph $G$ is defined by $Q^D(G)=Tr(G)+D(G)$ and $L^D(G)=Tr(G)-D(G)$, respectively.

In \cite{Plavsi}, Plavsi\'{c} et al. introduced the conception of \emph{Harary matrix} $RD(G)=(RD_{ij})$ of a graph $G$, which is also the so-called \emph{reciprocal distance matrix} and defined by
$$RD_{ij}=\left\{ {\begin{array}{*{20}{c}}
\frac{1}{d_{ij}},\;\text{if}\;\;i\neq j,\\
\;\;0,\;\text{otherwise}.
\end{array}}\right.$$
The \emph{reciprocal transmission} $RTr_{G}(v_i)$ of a vertex $v_i$ in $G$ is given by the sum of the reciprocal distance of $v_i$ from all other vertices of $G$, equivalently, $\sum_{v_j\in V(G)\setminus \{v_i\}}\frac{1}{d_{ij}}$. We say that a graph $G$ is \textit{$r$-reciprocal transmission regular} whenever its reciprocal transmission $RTr_{G}(v_i)=k$ for any $v_i\in V(G)$. Let $n\times n$ diagonal matrix $RT(G)$ be the vertex reciprocal transmissions of graph $G$. The \emph{reciprocal distance Laplacian matrix}, denoted by $RL(G)=RT(G)-RD(G)$, was defined by Bapat and Panda in \cite{Bapat}. It was proved \cite{Bapat} that, given a connected graph $G$ of order $n$, the spectral radius of its reciprocal distance Laplacian matrix  $\rho(RL(G))\leq n$ if and only if its complement graph, denoted by $\overline{G}$, is disconnected. The \emph{reciprocal distance signless Laplacian matrix}, denoted by $RQ(G)=RT(G)+RD(G)$, was introduced in \cite{Alhevaz}. Recently, the lower and upper bounds of the spectral radii of the reciprocal distance matrices and reciprocal distance signless Laplacian matrices of graphs were given in \cite{Das,Zhou} and \cite{Alhevaz,Medina}, respectively. Su et al. \cite{Su} determined the extremal graphs with maximal spectral radius of the $RD$-matrix among all graphs of order $n$ and precise vertex connectivity, edge connectivity, chromatic number and independence number, respectively. In \cite{Huang}, Huang et al. determined the extremal graphs with maximal spectral radius of the $RD$-matrix of the graphs with given matching number, bipartite graphs with given matching number, graphs with given cut edge number and so on.

In \cite{Nikiforov}, the matrix $A_{\alpha}(G)=\alpha \overline{D}(G)+(1-\alpha)A(G)$ for $\alpha\in [0,1]$, which is a convex linear combination of the degree diagonal matrix $\overline{D}(G)$ and adjacency matrix $A(G)$, was proposed by Nikiforov.
 Clearly, $A_{0}(G)=A(G)$, $A_{1/2}(G)={\frac{1}{2}}Q(G)$ and $A_{1}(G)=\overline{D}(G)$. In 2019, to study the gradual change from $D(G)$ to $Q^D(G)$, Cui et al. \cite{Cui2019} proposed  the matrix $D_{\alpha}(G)=\alpha Tr(G)+(1-\alpha)D(G)$ with $\alpha\in [0,1]$. It is easy to see that $D_{0}(G)=D(G)$, $D_{1/2}(G)=\frac{1}{2}Q^D(G)$ and $D_{1}(G)=Tr(G)$. Thus, the matrix $A_{\alpha}(G)$ is the unified way of $A(G)$ and $Q(G)$. In the same way, the matrix $D_{\alpha}(G)$ connects $D(G)$ to $Tr(G)$ with $\frac{1}{2}Q^D(G)$ situating in the middle of the interval. Recently, some spectral properties of $D_{\alpha}$-matrix  and $A_{\alpha}$-matrix have been studied extensively and several results have been published. For more details on $A_{\alpha}$-matrix and $D_{\alpha}$-matrix, see \cite{Cui2019,Diaz,Lin,Nikiforov} and the cited references therein.

In this paper, using a strategy similar to the one above,  the convex linear combinations of the matrices $RT(G)$ and $RD(G)$ are studied, which is defined by
$$RD_{\alpha}(G)=\alpha RT(G)+(1-\alpha)RD(G),\;\;0\leq \alpha\leq 1.$$
Since $RD_{0}(G)=RD(G)$, $RD_{1/2}(G)=\frac{1}{2}RQ(G)$ and $RD_{1}(G)=RT(G)$, then $RD_{1/2}(G)$ and $RQ(G)$ have same spectral properties. Thus $RD_{\alpha}(G)$ may form a unified theory of $RD(G)$ and $RQ(G)$. To this extent these matrices $RD(G)$, $RT(G)$ and $RQ(G)$ may be understood from a completely new perspective, and some interesting topics arise. Especially in spectral extremal graph theory, characterize the extremal graphs with maximal (or minimal) spectral radius among a given class of graphs is of great interest and importance. For the reciprocal distance matrix $RD(G)$, some spectral extremal graphs with fixed structure parameters have been characterized in \cite{Huang,Su}. It is natural to ask whether these result can be generalized to $RD_{\alpha}(G)$. If it is available, then we shall go straight to the extremal graphs with maximal (or minimal) spectral radius of  $RQ(G)$ for a graph $G$. Otherwise we can see some interesting differences.

Given a connected graph $G$, all eigenvalues of $RD_{\alpha}(G)$ forms the $RD_{\alpha}$-spectrum of $G$, use the notation $\sigma(RD_{\alpha}(G))=\{\lambda_{1}(RD_{\alpha}(G)),\lambda_{2}(RD_{\alpha}(G)),\ldots, \lambda_{n}(RD_{\alpha}(G))\}$, where all eigenvalues are arranged in descending order.
The maximum eigenvalue $\lambda_{1}(RD_{\alpha}(G))$ is called the spectral radius of the matrix $RD_{\alpha}(G)$, denoted by $\rho(RD_{\alpha}(G))$. Similarly, the spectrum of any matrix $M$ of order $n$ is denoted by $\sigma(M)=\{\lambda_{1}(M),\lambda_{2}(M),\ldots,\lambda_{n}(M)\}$ and its spectral radius is denoted by $\rho(M)$. In general, the notations $C_{n}$, $P_{n}$, $K_{n}$ and $K_{1,n-1}$ represent the cycle, path, complete graph and star of order $n$, respectively.

This paper is organized as follows. Section 2 presents some basics and spectra of $RD_{\alpha}(G)$ for some special graphs, such as join graphs, complete split graphs, complete multipartite graphs and so on. Section 3 presents some spectral properties of the matrix $RD_{\alpha}(G)$ and obtains some lower and upper bounds on spectral radius of $RD_{\alpha}(G)$ for a connected graph $G$. Section 4 is dedicated to spectral extremal problems, we determine the extremal graphs with maximal spectral radius of $RD_{\alpha}(G)$ in several types of simple connected graphs of order $n$ and precise vertex connectivity, edge connectivity, chromatic number and independence number, respectively. Finally, we sum up our previous work and put forward some
problems for further research in Section 5.

\section{Spectra of $RD_{\alpha}(G)$ for some special graphs}

In this section, we describe completely the $RD_{\alpha}$-spectra of some special graphs. At the first, we consider some easy cases. Clearly, $RD_{\alpha}(K_{n})=D_{\alpha}(K_{n})$. This implies that $\sigma(RD_{\alpha}(K_{n}))=\sigma(D_{\alpha}(K_{n}))=\{n-1,(\alpha n-1)^{[n-1]}\}$, where $a^{[b]}$ means that the multiplicity of $a$ is  $b$.

Given an $r$-regular graph $G$ of order $n$ with diameter $2$, based on the relationships between elements of $RD(G)$ and $A(G)$, we have $RD(G)=\frac{1}{2}(J_{n}-I_{n}+A(G))$, where $J_{n}$ and $I_{n}$ are the all ones matrix and identity matrix of size $n$, respectively. It is easy to get that $$RD_{\alpha}(G)=\frac{1}{2}[(\alpha n+\alpha r-1)I_{n}+(1-\alpha)J_{n}+(1-\alpha)A(G)].$$
Assume that $\{\lambda_1=r,\lambda_2,\ldots,\lambda_n\}$ is the spectrum of $A(G)$ and $\mathbf{1}_{n}=(1,1,\ldots,1)^{\rm T}$, $\textbf{x}_{2},\ldots,\textbf{x}_{n}$ are the corresponding eigenvectors. Since
\begin{equation*}
 RD_{\alpha}(G)\mathbf{1}_{n}=\frac{1}{2}(n+r-1)\mathbf{1}_{n}.
\end{equation*}
Then $\mathbf{1}_{n}=(1,1,\ldots,1)^{\rm T}$ is the eigenvector of $RD_{\alpha}(G)$ with respect to the eigenvalue $\frac{1}{2}(n+r-1)$. Noting that $\textbf{x}_{i}\perp\mathbf{1}_{n}$, then one has
\begin{equation*}
  \begin {aligned}
  RD_{\alpha}(G)\textbf{x}_{i}
  &=\frac{1}{2}((\alpha n+\alpha r-1)+(1-\alpha)\lambda_{i})\textbf{x}_{i}. \\
  \end {aligned}
\end{equation*}
Thus, $\frac{1}{2}((\alpha n+\alpha r-1)+(1-\alpha)\lambda_{i})$ is exactly an eigenvalue of $RD_{\alpha}(G)$ for arbitrary $i=2,3,\ldots,n$. So we can find the relationship between the $RD_{\alpha}$-eigenvalues of a regular graph with diameter $2$ and its adjacency eigenvalues.

In the following some known conclusions about a square matrix are given. We use $\Phi_{M}(\lambda)=\det(\lambda I_{n}-M)$ to represent the characteristic polynomial of $n\times n$ matrix $M$, which is also called as the $\lambda$\emph{-characteristic polynomial} of $M$. The \emph{coronal}, denoted by $\Gamma_{M}(\lambda)$, of a matrix $M$ is the sum of all the elements of the matrix $(\lambda I_{n}-M)^{-1}$. It means that
$$\Gamma_{M}(\lambda)=\mathbf{1}_{n}^{{\rm T}}(\lambda I_{n}-M)^{-1}\mathbf{1}_{n}.$$
It was proved \cite{Cui2012} that, if all the row sum of  matrix $M$ are equal  to $t$, then $ \Gamma_{M}(\lambda)=\frac{n}{\lambda-t}$. For two given column vectors $\textbf{x}$ and $\textbf{y}$, if a matrix $M$ is invertible, then $\det(M+\textbf{x}\textbf{y}^{{\rm T}})=(1+\textbf{y}^{{\rm T}}M^{-1}\textbf{x})\det(M)$, which is also called as the known \emph{Matrix Determinant Lemma} in \cite{Minc}.

Given two graphs $G_{1}=(V(G_{1}),E(G_{1}))$ and $G_{2}=(V(G_{2}),E(G_{2}))$, the graph $G_{1}\cup G_{2}$ has vertex set $V(G_{1}\cup G_{2})=V(G_{1})\cup V(G_{2})$ and edge set $E(G_{1}\cup G_{2})=E(G_{1})\cup E(G_{2})$. The \emph{join graph} $G_{1}\vee G_{2}$ is obtained from $G_{1}\cup G_{2}$ by connecting every vertex in $G_{1}$ to every vertex in $G_{2}$. In what follows, we give the $RD_{\alpha}$-characteristic polynomial of $G_{1}\vee G_{2}$.

\paragraph{Lemma 2.1.} For $i=1,2$, assume that $G_{i}$ is a connected graph of order $n_i$ and
$$M_{i}=\frac{1}{2}[(n_{i}\alpha-1)I_{n_{i}}+(1-\alpha)J_{n_{i}}+A_{\alpha}(G_i)],$$
where $A_\alpha(G_i)$ is the $A_\alpha$-matrix of $G_i$. Then the $RD_{\alpha}$-characteristic polynomial of $G_{1}\vee G_{2}$ is
\begin{equation}\label{0}
\Phi_{RD_{\alpha}}(\lambda)=\Phi_{M_{1}}(\lambda-\alpha n_{2})\Phi_{M_{2}}(\lambda-\alpha n_{1})[1-(1-\alpha)^{2}\Gamma_{M_{1}}(\lambda-\alpha n_{2})\Gamma_{M_{2}}(\lambda-\alpha n_{1})].
\end{equation}
\begin{proof}
Since the diameter of $G_{1}\vee G_{2}$ is $2$, then
$$RD_{ij}(G_{1}\vee G_{2})=
\left\{ \begin{array}{cc}
          1 & \text{if}\;\; v_{i}v_{j}\in E(G_{1}\vee G_{2}), \\
          0 & \text{if}\;\; i=j, \\
          \frac{1}{2} & \text{otherwise}.
        \end{array}
\right.$$
With a suitable label of the vertices of $G_{1}\vee G_{2}$, one gets
$$RD_{\alpha}(G_{1}\vee G_
{2})=
\left(
    \begin{array}{cc}
      \alpha n_{2}I_{n_{1}}+M_{1} & (1-\alpha) J_{n_{1}\times n_{2}}\\
      (1-\alpha) J_{n_{2}\times n_{1}} &  \alpha n_{1}I_{n_{2}}+M_{2}
    \end{array}
  \right).$$
The remaining proof is exactly similar to that of Lemma 2.1 in \cite{Lin}, we omit the detail.
\end{proof}

\paragraph{Theorem 2.2.} For $i=1,2$, suppose that $G_{i}$ is an $r_{i}$-regular graph with $n_i$ vertices. Then the spectrum of $RD_{\alpha}(G_{1}\vee G_{2})$ consists precisely of:

\begin{enumerate}[(i)]
  \item The eigenvalue $\frac{1}{2}[\alpha (n + {n_2} + {r_1}) + (1 - \alpha ){\lambda _j}(A({G_1})) - 1] \;\;\text{for any} \;\;2\leq j\leq n_{1}$;
  \item The eigenvalue $\frac{1}{2}[\alpha (n + {n_1} + {r_2}) + (1 - \alpha ){\lambda _k}(A({G_2})) - 1] \;\; \text{for any} \;\;2\leq k\leq n_{2}$;
  \item The rest eigenvalues are the roots of  $(\lambda-\alpha n_{2}-\frac{1}{2} n_{1}-\frac{1}{2}r_{1}+\frac{1}{2})(\lambda-\alpha n_{1}-\frac{1}{2} n_{2}-\frac{1}{2}r_{2}+\frac{1}{2})-(1-\alpha)^{2}n_{1}n_{2}=0$.
\end{enumerate}

\begin{proof}
For $i=1,2$, we take $M_{i}=\frac{1}{2}[(n_{i}\alpha-1)I_{n_{i}}+(1-\alpha)J_{n_{i}}+A_{\alpha}(G_i)]$, where $A_{\alpha}(G_i)$ is the $A_\alpha$-matrix of $G_i$. First we calculate
${\Phi_{{M_1}}}(\lambda-\alpha n_2 ) = \det((\lambda-\alpha n_2 ) {I_{{n_1}}} - \frac{1}{2}(({n_1}\alpha  - 1){I_{{n_1}}} + (1 - \alpha ){J_{{n_1}}} + {A_\alpha }(G_1)))$. For the sake of convenience, let ${N_1} = \frac{1}{2}(\alpha{n_1}+2 \alpha{n_2} - 1){I_{{n_1}}} + \frac{1}{2}{A_\alpha }(G_1)$. Now, from the Matrix Determinant Lemma, we can obtain
\begin{equation}\label{1}
  \begin {aligned}
{\Phi_{{M_1}}}(\lambda-\alpha n_2) &= (1 - \frac{1}{2}(1 - \alpha ){\textbf{1}_{n_1}^{\rm T}}{(\lambda {I_{{n_1}}} - {N_1})^{ - 1}}{\textbf{1}_{{n_1}}})\det (\lambda {I_{{n_1}}} - {N_1})\\
  &= (1 - \frac{1}{2}(1 - \alpha ){\Gamma _{{N_1}}}(\lambda )){(\frac{1}{2})^{{n_1}}}{\Phi_{{A_\alpha }({G_1})}}(2\lambda  - \alpha {n_1} - 2\alpha {n_2} + 1).
  \end {aligned}
\end{equation}
Since $G_{1}$ is an $r_{1}$-regular graph, then the matrix $N_{1}$ has the same row sum $\frac{1}{2}(\alpha n_{1}+2\alpha {n_2}+r_{1}-1)$. This implies that
\begin{equation}\label{2}
{\Gamma _{{N_1}}}(\lambda )={\textbf{1}_{n_1}^{\rm T}}{(\lambda {I_{{n_1}}} - {N_1})^{ - 1}}{\textbf{1}_{{n_1}}}= \frac{{{n_1}}}{{\lambda  - \frac{1}{2}(\alpha {n_1}+2\alpha {n_2} + {r_1} - 1)}}.
\end{equation}
Substituting (\ref{2}) into (\ref{1}), one has
\begin{equation}\label{3}
{\Phi_{{M_1}}}(\lambda  - \alpha {n_2}) = {(\frac{1}{2})^{{n_1}}}\frac{{\lambda  -\alpha {n_2}- \frac{1}{2}{n_1} - \frac{1}{2}{r_1} + \frac{1}{2}}}{{\lambda  - \frac{1}{2}(\alpha {n_1} +2\alpha {n_2}+ {r_1} - 1)}}{\Phi_{{A_\alpha }({G_1})}}(2\lambda  - \alpha {n_1}-2\alpha {n_2} + 1).
\end{equation}
Similarly,
\begin{equation}\label{4}
{\Phi_{{M_2}}}(\lambda  - \alpha {n_1}) = {(\frac{1}{2})^{{n_2}}}\frac{{\lambda  -\alpha {n_1}- \frac{1}{2}{n_2} - \frac{1}{2}{r_2} + \frac{1}{2}}}{{\lambda  - \frac{1}{2}(\alpha {n_2} +2\alpha {n_1}+ {r_2} - 1)}}{\Phi_{{A_\alpha }({G_2})}}(2\lambda  - \alpha {n_2}-2\alpha {n_1} + 1).
\end{equation}
Since $G_{i}$ is an $r_{i}$-regular graph for $i=1,2$. Then the matrix $M_{i}$ has the same row sum $\frac{1}{2}(n_{i}+r_{i}-1)$. Thus,
\begin{equation}\label{5}
{\Gamma _{{M_1}}}(\lambda-\alpha n_2 ) = \frac{{{n_1}}}{{\lambda -\alpha n_2 - \frac{1}{2}({n_1} + {r_1} - 1)}},\;\;{\Gamma _{{M_2}}}(\lambda-\alpha n_1 ) = \frac{{{n_2}}}{{\lambda -\alpha n_1 - \frac{1}{2}({n_2} + {r_2} - 1)}}.
\end{equation}
Substituting (\ref{3}), (\ref{4}), (\ref{5}) back into (\ref{0}) in Lemma 2.1, we obtain the $RD_{\alpha}$-characteristic polynomial of $G_{1}\vee G_{2}$
\[
\Phi_{RD_{\alpha}}(\lambda)=(\frac{1}{2})^{n_{1}+n_{2}}\frac{\Phi_{A_{\alpha}(G_{1})}(2\lambda-2\alpha n_{2}-\alpha n_{1}+1)\Phi_{A_{\alpha}(G_{2})}(2\lambda-2\alpha n_{1}-\alpha n_{2}+1)}{(\lambda-\alpha n_{2}-\frac{1}{2}\alpha n_{1}-\frac{1}{2}r_{1}+\frac{1}{2})(\lambda-\alpha n_{1}-\frac{1}{2}\alpha n_{2}-\frac{1}{2}r_{2}+\frac{1}{2})}f(\lambda),
\]
where $f(\lambda) = (\lambda  - \alpha {n_2} - \frac{1}{2}{n_1} - \frac{1}{2}{r_1} + \frac{1}{2})(\lambda  - \alpha {n_1} - \frac{1}{2}{n_2} - \frac{1}{2}{r_2} + \frac{1}{2}) - {(1 - \alpha )^2}{n_1}{n_2}$.
Notice that $A_{\alpha}(G_{i})=\alpha r_{i}I_{n_{i}}+(1-\alpha)A(G_{i})$ for the $r_i$-regular graph $G_{i}$. Thus
$\alpha r_{i}+(1-\alpha)\lambda$ is an eigenvalue of $A_{\alpha}(G_{i})$ whenever $\lambda$ is an eigenvalue of $G_{i}$  for  $i=1,2$. Therefore, the required result follows.
\end{proof}

Obviously, Theorem 2.2 implies immediately the $RD_{\alpha}$-spectra of the following graphs.

\paragraph{Corollary 2.3.} For $a,b\geq 1$, let $K_{a,b}=\overline{K_{a}}\vee\overline{K_{b}}$ be a complete bipartite graph with $a+b=n$ vertices. Then the spectrum of $RD_{\alpha}(K_{a,b})$ consists precisely of:
\[
(\frac{\alpha(n+b)-1}{2})^{[a-1]}, \;\;\;\;\;\;(\frac{\alpha(n+a)-1}{2})^{[b-1]},
\]
and the rest eigenvalues are
\[
\frac{{(\alpha  + \frac{1}{2})n - 1 \pm \sqrt {{{(\alpha  - \frac{1}{2})}^2}{{(a - b)}^2} + 4{{(1 - \alpha )}^2}ab} }}{2}.
\]

\paragraph{Corollary 2.4.} For $a\geq1, b\geq 2$, let $CS_{a,b}=K_{a}\vee\overline{K_{b}}$ be a complete split graph with $a+b=n$ vertices. Then the spectrum of $RD_{\alpha}(CS_{a,b})$ is given by
$(\alpha n - 1)^{[a-1]}$, $(\frac{{\alpha (n + a) - 1}}{2})^{[b-1]}$ and
\[
\frac{{(\alpha  + 1)n - \frac{1}{2}b - \frac{3}{2} \pm \sqrt {{{((\alpha  - 1)(a - b) - \frac{1}{2}b + \frac{1}{2})}^2} + 4{{(1 - \alpha )}^2}ab} }}{2}.
\]

A wheel graph $W(n)$ is the graph constructed by the join operation between a isolated vertex $K_{1}$ and a cycle $C_{n-1}$, that is, $W(n)=K_{1}\vee C_{n-1}$.

\paragraph{Corollary 2.5.} The $RD_{\alpha}$-eigenvalues of the wheel graph $W(n)=K_{1}\vee C_{n-1}$ are
\[
\frac{{\alpha (n + 3) - 1 + 2(1 - \alpha )\cos (\frac{{2j\pi }}{{n - 1}})}}{2},\;\; \text{for} \;\;1\leq j\leq n-2,
\]
and
\[
\frac{{(\alpha  + \frac{1}{2})n \pm \sqrt {{{(\alpha (n - 2) - \frac{n}{2})}^2} + 4{{(1 - \alpha )}^2}(n - 1)} }}{2}.
\]

The following proposition establishes a relationship between the  $RD_\alpha$-spectrum and some particular vertex subset of $G$, which will is used to give the spectrum of a graph with given clusters. The proof of this proposition is similar to that of Proposition 11 in \cite{Cui2019}, the detail is omitted.

\paragraph{Proposition 2.6.} Suppose that $G$ be a connected graph and $C\subseteq V(G)$ such that any two vertices in $C$ having same neighborhood in $V(G)\setminus C$, where $|C|=c$. Then
\begin{enumerate}[(i)]
  \item If $C$ is an independent set, then $\alpha t+\frac{\alpha-1}{2}$ is an eigenvalue of $RD_{\alpha}(G)$ with repeated at least $c-1$ times, where $t$ is the reciprocal transmission in $G$ of the vertices in $C$.
  \item If $C$ is a clique, then $\alpha(t+\frac{1}{2}c+\frac{1}{2})-1$ is an eigenvalue of $RD_{\alpha}(G)$ with repeated at least $c-1$ times.
\end{enumerate}

The following definition comes from \cite{Cardoso2017,Merris}. Let $G$ be a connected graph of order $n$. A \emph{cluster} in graph $G$, denoted by $(C,S)$,  is a pair of vertex subsets $C$ and $S$, where  $|C| =c$, $|S| =s$, and $C$ is a set of cardinality $|C| =c \geq2$ of pairwise co-neighbor vertices sharing the same set $S $ of $s$ neighbors. Clearly, $C$ is an independent set of $G$ and the reciprocal transmissions $RTr_G(v)$ are the same for any $v\in C$. Graph $G(K_c)$ is the graph constructed by  replacing the independent set $C$ in graph $G$ with the complete graph $K_c$.

For simplicity, let $C=\{v_1,\ldots,v_c\}$, $S=\{v_{c+1},\ldots,v_{c+s}\}$ and $\{v_{c+s+1},\ldots,v_{n}\}$ be the set of the remaining vertices in $G$. Also, let $t$ be the reciprocal transmission of the vertices in $C$ and $\beta=1-\alpha$. Then the matrices $RD_{\alpha}(G)$ and $RD_{\alpha}(G(K_{c}))$ can be given by:
\[
R{D_\alpha }(G) = \left( {\begin{array}{*{20}{c}}
  U&W \\
  {{W^{\rm T}}}&Z
\end{array}} \right)
\;\;\text{and}\;\;
R{D_\alpha }(G({K_{c}})) = \left( {\begin{array}{*{20}{c}}
  V&W \\
  {{W^{\rm T}}}&Z
\end{array}} \right),
\]
where
\[
U=U_{c\times c} = \left( {\begin{array}{*{20}{c}}
  {\alpha t}&{\frac{1}{2}\beta }& \ldots &{\frac{1}{2}\beta } \\
  {\frac{1}{2}\beta }&{\alpha t }& \ddots & \vdots  \\
   \vdots & \ddots & \ddots &{\frac{1}{2}\beta } \\
  {\frac{1}{2}\beta }& \ldots &{\frac{1}{2}\beta }&{\alpha t}
\end{array}} \right),
\]
\[
V=V_{c\times c} = \left( {\begin{array}{*{20}{c}}
  {\alpha (t  + \frac{1}{2}c - \frac{1}{2})}&\beta & \ldots &\beta  \\
  \beta &{\alpha (t  + \frac{1}{2}c - \frac{1}{2})}& \ddots & \vdots  \\
   \vdots & \ddots & \ddots &\beta  \\
  \beta & \ldots &\beta &{\alpha (t + \frac{1}{2}c - \frac{1}{2})}
\end{array}} \right),
\]
\[
W = \beta \left( {\begin{array}{*{20}{c}}
  {\frac{1}{{{d_{1,c + 1}}}}}&{\frac{1}{{{d_{1,c + 2}}}}}& \ldots &{\frac{1}{{{d_{1,n}}}}} \\
  {\frac{1}{{{d_{1,c + 1}}}}}&{\frac{1}{{{d_{1,c + 2}}}}}& \ldots &{\frac{1}{{{d_{1,n}}}}} \\
   \vdots & \vdots & \vdots & \vdots  \\
  {\frac{1}{{{d_{1,c + 1}}}}}&{\frac{1}{{{d_{1,c+ 2}}}}}& \ldots &{\frac{1}{{{d_{1,n}}}}}
\end{array}} \right)
\]
and
\[
Z=Z_{(n-c)\times(n-c)} = \left( {\begin{array}{*{20}{c}}
  {\alpha RT_G(v_{c + 1})}&{\frac{\beta }{{{d_{c + 1,c + 2}}}}}& \ldots &{\frac{\beta }{{{d_{c + 1,n}}}}} \\
  {\frac{\beta }{{{d_{c + 2,c + 1}}}}}&{\alpha RT_G(v_{c + 2})}& \ldots &{\frac{\beta }{{{d_{c + 2,n}}}}} \\
   \vdots & \vdots & \ddots & \vdots  \\
  {\frac{\beta }{{{d_{n,c + 1}}}}}&{\frac{\beta }{{{d_{n,c + 2}}}}}& \ldots &{\alpha RT_G(v_{n})}
\end{array}} \right).
\]

\paragraph{Theorem 2.7.} Let $G$ be a connected graph with $n$ vertices and having a cluster $(C,S)$. Assume that the reciprocal transmissions $RTr_{G}(v)=t$ for any $v\in C$. Then we have the following statement.
\begin{enumerate}[(i)]
  \item $\alpha t+\frac{\alpha-1}{2}$ is an eigenvalue of $RD_{\alpha}(G)$ with multiplicity $c-1$, and the rest $n -c + 1$  eigenvalues of $RD_{\alpha}(G)$ are the eigenvalues of the following matrix
  \[
  Q_1 = \left( {\begin{array}{*{20}{c}}
	{\alpha (t + \frac{1}{2}) - \frac{1}{2} + \frac{1}{2}c ( {1 - \alpha } )}&{{\textbf{x}^{\rm T}}} \\
	c\textbf{x}&Z
	\end{array}} \right),
  \]
  where
  \[
  \textbf{x} = (1 - \alpha )\left(
  {\frac{1}{{{d_{1,c + 1}}}}},\;{\frac{1}{{{d_{1,c + 2}}}}},\;\ldots,\;{\frac{1}{{{d_{1,n}}}}}\right)^{\rm T}.
  \]
  \item $\alpha(t+\frac{1}{2}c+\frac{1}{2})-1$ is an eigenvalue of $RD_{\alpha}(G(K_{c}))$ with repeated $c-1$ times, and the rest $n -c + 1$ eigenvalues of $RD_{\alpha}(G(K_c))$ are the eigenvalues of the following  matrix
  \[
  {Q_2} = \left( {\begin{array}{*{20}{c}}
  	{\alpha (t + \frac{1}{2}) - 1 +\frac{1}{2} c ( {2 - \alpha } )}&{{\textbf{x}^{\rm T}}} \\
  	{c\textbf{x}}&Z
  	\end{array}} \right).
  \]
\end{enumerate}

\begin{proof} Proposition 2.6 implies that $\alpha t+\frac{\alpha-1}{2}$ is an eigenvalue of $RD_{\alpha}(G)$ with repeated $c-1$ times. Similarly, $\alpha(t+\frac{1}{2}c+\frac{1}{2})-1$ is an eigenvalue of $RD_{\alpha}(G(K_{c}))$ with repeated $c-1$ times.
	
In what follows, suppose that $\pi: V(G)={V_1}\cup\{v_{c+1}\}\cup\{v_{c+2}\}\cup\cdots\cup\{v_n\}$ is a partition of the vertix set $V(G)$, where the vertex subset ${V_1} =C= \{v_1,\ldots,v_c\}$. By observation, this partition $\pi$ is an equitable partition for the matrix $RD_{\alpha}(G)$ and $RD_{\alpha}(G(K_{c}))$. It is easy to see that $Q_1$ and $Q_2$ are the quotient matrices corresponding to $RD_{\alpha}(G)$ and $RD_{\alpha}(G(K_{c}))$ with respect to the partition $\pi$. Therefore, the eigenvalues of $Q_1$ and $Q_2$ are the $( {n - c + 1})$ remaining eigenvalues of $RD_{\alpha}(G)$ and $RD_{\alpha}(G(K_{c}))$, respectively.
\end{proof}

A vertex with only one neighbor is called a \textit{pendent vertex} in a graph $G$. The unique neighbor of a pendent vertex is called a \textit{quasi-pendent vertex} in $G$. Clearly, a pair of each maximal set $C$ of pendent vertices and the corresponding quasi-pendent vertex $v$ forms a cluster $(C,\{v\})$.

\paragraph{Remark 1.} Taking $\alpha  = 0$ and $\alpha  = \frac{1}{2}$ in Theorem 2.7, we can get some spectral information about the reciprocal distance matrix $RD(G)$ and reciprocal distance signless Laplacian of a graph with clusters, respectively. For example, let graph $G$  with a cluster $(C,S)$ is a connected graph and cardinality $|C|=c$, Theorem 2.7 implies that $-\frac{1}{2}$ is an eigenvalue of $RD(G)$ repeated at least $c-1$ times. Now assume that $G$ has $p$ pendent vertices and $q$ quasi-pendent vertices. Apply Theorem 2.7 repeatedly, it turns out that $-\frac{1}{2}$ is an eigenvalue of $RD(G)$ with repeated at least $p - q$ times.\\

Next we give the spectrum of $RD_{\alpha}(G)$ for a complete $r$-partite graph $G=K_{n_1,n_2,\ldots,n_r}$ with $n=n_1+n_2+\cdots+n_r$ vertices. A complete $r$-partite graph $G=K_{n_1,n_2,\ldots,n_r}$ is the graph formed by partitioning its vertex set into $r$ subsets $V_1,V_2,\ldots,V_r$ with cardinality $|V_i|=n_i$, and joining two vertices by an edge if and only if they are located in different subsets.

\paragraph{Theorem 2.8.}  For a complete $r$-partite graph $G=K_{n_1,n_2,\ldots,n_r}$ of order $n=n_1+n_2+\cdots+n_r$ with $n \geqslant 4$, we have
\begin{enumerate}[(i)]
  \item $\alpha (n - \frac{n_i}{2}) - \frac{1}{2}$ is an eigenvalue of $RD_{\alpha}(G)$ with repeated $n_i-1$ times, for any $i\in \left\{ {1,2, \ldots ,r} \right\}$.
  \item the eigenvalues of the following matrix $T$ are the rest  eigenvalues of $RD_{\alpha}(G)$, where
\[
T = \left( {\begin{array}{*{20}{c}}
  {\alpha (n - n_1) + \frac{1}{2}(n_1 - 1)}&{(1 - \alpha )n_2}& \ldots &{(1 - \alpha )n_r} \\
  {(1 - \alpha )n_1}&{\alpha (n - n_2) + \frac{1}{2}(n_2- 1)}& \ldots &{(1 - \alpha )n_r} \\
   \vdots & \vdots & \ddots & \vdots  \\
  {(1 - \alpha )n_1}&{(1 - \alpha )n_2}& \ldots &{\alpha (n - n_r) + \frac{1}{2}(n_r- 1)}
\end{array}} \right).
\]
\end{enumerate}

\begin{proof} Let $\pi:{V_1} \cup {V_2} \cup  \ldots  \cup {V_r}$ be the partition of vertex set $V(G)$ with cardinality $|V_i|=n_i$ for $i\in \left\{ {1,2, \ldots ,r} \right\}$. Notice that $RTr_G(v)=n - \frac{n_i}{2}-\frac{1}{2}$ for any vetrex $v\in V_i$. It follows from Proposition 2.6 that $\alpha (n - \frac{n_i}{2}) - \frac{1}{2}$ is an eigenvalue of $RD_{\alpha}(G)$ with repeated $n_i-1$ times. On the other hand, the partition $\pi$ is an equitable partition of $RD_{\alpha}(G)$ for the complete $r$-partite graph $G$. With a suitable labeling, the quotient matrix of $RD_{\alpha}(G)$ can be described as the matrix $T$ in the satement of this theorem. Hence, the $r$ remaining eigenvalues of $RD_{\alpha}(G)$ can be derived from the eigenvalues of $T$.
\end{proof}

\section{Basic properties of $RD_{\alpha}(G)$}

\paragraph{Lemma 3.1.\cite{So}} For two Hermitian matrices $A_{n\times n}$ and $B_{n\times n}$, let $C=A+B$ and $1\leq i,j\leq n$. Then
\[
\begin{gathered}
  {\lambda _i}(C) \leq {\lambda _j}(A) + {\lambda _{i - j + 1}}(B)\;\; {\rm for}\;\; j \leq i, \hfill \\
  {\lambda _i}(C) \geq {\lambda _j}(A) + {\lambda _{i - j + n}}(B)\;\; {\rm for}\;\; i \leq j. \hfill \\
\end{gathered}
\]
Moreover, each one of the equality holds iff there is a nonzero vector that is an eigenvector to each of the three eigenvalues involved.

Lemma 3.1 implies that
\[
{\lambda _i}(A) + {\lambda _{\min }}(B) \leq {\lambda _i}(C) \leq {\lambda _i}(A) + {\lambda _{\max }}(B).
\]

The following proposition implies that, for $\frac{1}{2}\leq\alpha\leq 1$, each one of the eigenvalues of $RD_{\alpha}(G)$ does not decrease when an edge is added between two nonadjacent vertices in a graph $G$.

\paragraph{Proposition 3.2.} Suppose that $G$ is a graph with $n$ vertices and $\frac{1}{2}\leq\alpha\leq 1$. For an edge $e\notin E(G)$, let the graph $\widetilde{G}=G+e$. Then we have
\begin{equation}\label{7}
 {\lambda _i}(R{D_\alpha }(\widetilde G)) \geq {\lambda _i}(R{D_\alpha }(G))\;\;{\rm for}\;\; 1\leq i\leq n.
\end{equation}
\begin{proof}
For any two vertices $u$ and $v$ of graph $G$, we find that $d_{\widetilde{G}}(u,v)\leq d_{G}(u,v)$, equivalently, $\frac{1}{d_{\widetilde{G}}(u,v)}\geq \frac{1}{d_{G}(u,v)}$. Then $RD_{\alpha}(\widetilde{G})=RD_{\alpha}(G)+M$, where
\[
M = \left( {\begin{array}{*{20}{c}}
  {\alpha {m_1}}&{(1 - \alpha ){m_{1,2}}}& \cdots &{(1 - \alpha ){m_{1,n}}} \\
  {(1 - \alpha ){m_{2,1}}}&{\alpha {m_2}}& \cdots &{(1 - \alpha ){m_{2,n}}} \\
   \vdots & \vdots & \ddots & \vdots  \\
  {(1 - \alpha ){m_{n,1}}}&{(1 - \alpha ){m_{n,2}}}& \cdots &{\alpha {m_n}}
\end{array}} \right)
\]
and ${m_i} = \sum\nolimits_{j = 1,j \ne i}^{n} {{m_{i,j}}} $ for $i=1,2,\ldots, n$. It is clear that, for $\alpha\in [\frac{1}{2},1]$, the matrix $M$ is diagonally dominant. At that time, $M$ is also a positive semidefinite matrix with nonnegative diagonal entries. Thus the minimum eigenvalue of $M$ is greater than or equal to zero. Now, Lemma 3.1 implies that the required result follows.
\end{proof}

In the light of Proposition 3.2, an upper bound for the $k$-th largest eigenvalue of $RD_{\alpha}(G)$ for a connected graph $G$ of order $n$ is given by
\[
{\lambda _i}(R{D_\alpha }(G))\leq{\lambda _i}(R{D_\alpha }({K_n}))= \alpha n - 1.
\]
where $2\leq i\leq n$ and $\frac{1}{2}\leq\alpha\leq 1$.

What is noteworthy is that the inequality (\ref{7}) is strict for $i=1$ and $\alpha\in [0,1)$. Indeed, given a connected graph $G$,  $RD_{\alpha}(G)$ is a nonnegative irreducible matrix for $\alpha\in [0,1)$. If an edge is added between two nonadjacent vertices of $G$, then at least one entry of $RD_{\alpha}(G)$ strictly increases. Perron-Frobenius theory implies that the spectral radius of $RD_{\alpha}(G)$ strictly increases whenever we add an edge $e\notin E(G)$ to graph $G$. Thus we shall see next.

\paragraph{Proposition 3.3.}
Suppose that $G$ is a connected graph. For an edge $e\notin E(G)$, let $\widetilde{G}=G+e$. Then, for $\alpha\in [0,1)$,
\begin{equation}\label{8}
\rho (R{D_\alpha }(\tilde G)) > \rho (R{D_\alpha }(G)).
\end{equation}

\paragraph{Corollary 3.4.}Let $G$ be a connected bipartite graph with  $n$ vertices, where $n\geq2$. Suppose that the two disjoint parts have $a$ and $n-a$ vertices with $1\leq a\leq \lfloor\frac{n}{2}\rfloor$, then
\[
\rho (R{D_\alpha }(G))\leq \frac{{(\alpha  + \frac{1}{2})n - 1 + \sqrt {{{(\alpha  - \frac{1}{2})}^2}{{(2a - n)}^2} + 4{{(1 - \alpha )}^2}a(n - a)} }}{2}.
\]
Then the equality holds if and only if $G=K_{a,n-a}$.

\begin{proof}
It follows from Proposition 3.3 that $\rho (R{D_\alpha }(K_{a,n-a}) \geq \rho (R{D_\alpha }(G))$. Moreover, the equality holds if and only if $G=K_{a,n-a}$. On the basis of Corollary 2.3, we have
\[
\rho (R{D_\alpha }(K_{a,n-a})=\frac{{(\alpha  + \frac{1}{2})n - 1 + \sqrt {{{(\alpha  - \frac{1}{2})}^2}{{(2a - n)}^2} + 4{{(1 - \alpha )}^2}a(n - a)} }}{2},
\]
as required.
\end{proof}

The following proposition 3.5 implies that, for any connected graph $G$, the eigenvalue $\lambda_{i}(RD_{\alpha}(G))$ is increasing in $\alpha$, where $0\leq\alpha\leq1$ and $2\leq i\leq n$.

\paragraph{Proposition 3.5.}
Let $G$ be a connected graph with $n$ vertices. If $1\geq \alpha > \beta \geq 0$ and  $1\leq i \leq n$, then
\begin{equation}\label{9}
{\lambda _i}(R{D_\alpha }(G)) \geq {\lambda _i}(R{D_\beta }(G)).
\end{equation}
The equality holds in (\ref{9}) if and only if $i=1$ and graph $G$ is  reciprocal transmission regular.
\begin{proof}
As we all known, $RD_\alpha(G) = RD_\beta(G) + (\alpha  - \beta )(RT(G) - RD(G))=RD_\beta(G) + (\alpha  - \beta )(RL(G))$.
Applying Lemma 3.1 with $C = RD_\alpha(G)$, $A = RD_\beta(G)$ and $B = (\alpha-\beta)RL(G)$, we have
\[
{\lambda _i}(R{D_\alpha }(G)) \geq {\lambda _i}(R{D_\beta }(G)) +(\alpha  - \beta ) \lambda_{n}(RL(G)).
\]
Since the matrix $RL(G)$ is positive semidefinite, then ${\lambda _n}(RL(G))\geq 0$ (in fact, ${\lambda _n}(RL(G))= 0$ \cite{Bapat}). Hence,
\[
{\lambda _i}(R{D_\alpha }(G)) \geq {\lambda _i}(R{D_\beta }(G)).
\]

Now suppose that the equality holds in (\ref{9}). According to Lemma 3.1, ${\lambda _i}(R{D_\alpha }(G))$, ${\lambda _i}(R{D_\beta }(G))$ and $\lambda_{n} (RL(G))$ have the same eigenvector $\textbf{1}_{n}$ because $\textbf{1}_{n}$ is the unique eigenvector of ${\lambda _n}(RL(G))$ (see Remark 2.1 in \cite{Bapat}). This means that $i=1$ and the graph $G$ is reciprocal transmission regular. To the contrary, let $i=1$ and $G$ be reciprocal transmission regular. Then, it can be verified that the equality in (\ref{9}) holds.
\end{proof}

Proposition 3.5 implies that the following corollary is immediate.

\paragraph{Corollary 3.6.} Let $G$ be a connected graph with $n$ vertices and $0\leq\alpha\leq1$. Then
\[
{\lambda _k}(RD(G)) \leq {\lambda _k}(R{D_\alpha }(G)) \leq {\lambda _k}(RT(G)) = RT{r_k} \;\;\text{for}\;\; 1\leq k \leq n,
\]
where $RT{r_k}$ is the $k$-th largest reciprocal transmission of $G$.

Remark that Corollary 3.6 gives an upper bound on $\lambda_{n}(RD_{\alpha}(G))$, but this bounds can be improved in the following proposition.

\paragraph{Proposition 3.7.} Suppose that $G$ is a connected graph with $n$ vertices. Then
\[
{\lambda _n}(R{D_\alpha }(G)) \leq \alpha RT{r_n} \;\; {\rm and} \;\;\rho (R{D_\alpha }(G)) \geq \alpha RT{r_1}
\]
where $RTr_{1}$ and $RTr_{n}$ are the maximum reciprocal transmission and  minimum reciprocal transmission of $G$, respectively.

\begin{proof}
Assume that $u$ is a vertex with the minimum reciprocal transmission $RTr_{n}$ of graph $G$. If $\textbf{e}_u$ is the characteristic vector of the vertex $u$, then
\[
{\lambda _n}(R{D_\alpha }(G)) = \mathop {\min }\limits_{{{\left\| \textbf{x} \right\|}_2} = 1} {\textbf{x}^{\rm T}}R{D_\alpha }(G)\textbf{x} \leq {\textbf{e}_u^{\rm T}}R{D_\alpha }(G)\textbf{e}_u=\alpha RTr_{n}.
\]
Similarly, we easily obtain that $\rho (R{D_\alpha }(G)) \geq \alpha RT{r_1}$.
\end{proof}

Further, Using Lemma 3.1 to the definition of $R{D_\alpha }(G)$, we can get the following bounds on the $k$-th largest eigenvalue ${\lambda _k}(R{D_\alpha }(G))$ in terms of the $k$-th largest eigenvalue ${\lambda _k}(R{D }(G))$, maximum reciprocal transmission $RTr_{1}$ and minimum reciprocal transmission $RTr_{n}$ of a connected graph $G$. That is to say,
\[
\alpha RT{r_n} + (1 - \alpha ){\lambda _k}(RD(G)) \leq {\lambda _k}(R{D_\alpha }(G)) \leq \alpha RT{r_1} + (1 - \alpha ){\lambda _k}(RD(G)).
\]

Observe that $R{D_\alpha }(G) + R{D_{1 - \alpha }}(G) = RQ(G)$. From Lemma 3.1 again, we can obtain a relation for the spectral radii of three matrices below:
\begin{equation}\label{10}
\rho (R{D_\alpha }(G)) + \rho (R{D_{1 - \alpha }}(G)) \geq\rho (RQ(G)).
\end{equation}
Remark that this relation can been performed the transformation between the upper and lower bounds of the spectral radius $\rho(R{D_\alpha }(G))$. For example, applying this relation and Lemma 3.1, we can get the following inequalities on $\rho (R{D_\alpha }(G))$ in terms of $\rho (RD(G))$, $\rho (RQ(G))$ and the maximum reciprocal transmission $RTr_{1}$ of $G$. This statement and its proof is analogous to an existing result related
to $\rho (D_{\alpha}(G))$ (see Propositions 7 and 8 in \cite{Diaz}), these details of the proof are omitted.

\paragraph{Proposition 3.8.}
Let $G$ be a connected graph of order $n$ with maximum reciprocal transmission $RTr_{1}$. Then

\begin{enumerate}[(i)]
\item For $0\leq \alpha \leq \frac{1}{2}$, we have
\[
(1-\alpha) \rho (RQ(G)) + (2\alpha -1)RTr_{1}\leq
\rho (R{D_\alpha }(G)) \leq \alpha \rho (RQ(G)) + (1 - 2\alpha )\rho (RD(G)).
\]
Moreover, if $G$ is not reciprocal transmission regular, then the first equality holds if and only if $\alpha=\frac{1}{2}$, and the second equality holds if and only if $\alpha=0$ or $\alpha=\frac{1}{2}$.
\item For $\frac{1}{2} \leq \alpha \leq 1$, we have
\[\alpha\rho (RQ(G)) + ( 1-2\alpha )\rho (RD(G))\leq\rho (R{D_\alpha }(G)) \leq (1-\alpha) \rho (RQ(G)) + ( 2\alpha -1)RTr_{1}\]
Moreover, if $G$ is not reciprocal transmission regular, then the first equality holds if and only if $\alpha=\frac{1}{2}$, and the second equality holds if and only if $\alpha=\frac{1}{2}$ or $\alpha=1$.
\end{enumerate}

In what follows, we turn our attention to the positive semidefiniteness of the matrix $R{D_\alpha }(G)$. First notice that the matrix $R{D_0 }(G)$ is not positive semidefinite, but the matrix $R{D_\frac{1}{2}}(G)$ is positive semidefinite.
According to the Proposition 3.5, the function ${f}(\alpha )= {\lambda _{\min }}(R{D_\alpha }(G))$ is continuous and increasing in $\alpha$.
Then there exists the minimum $\alpha_{0}\in(0,\frac{1}{2}]$ such that ${f}(\alpha_{0} )=0$. Thus the matrix $R{D_\alpha }(G)$ is a positive semidefinite matrix if and only if $\alpha\in[\alpha_{0}, 1]$. This observation raises the following problem:

\paragraph{Problem 3.9.} For a connected graph, determine the minimum $\alpha_{0}\in(0,\frac{1}{2}]$ such that the matrix $R{D_\alpha }(G)$ is a positive semidefinite matrix for $\alpha\in[\alpha_{0}, 1]$.

\paragraph{Proposition 3.10.}
 Let $G$ be a $k$-reciprocal transmission regular graph with $n$ vertices. Then the matrix $R{D_\alpha }(G)$ is positive semidefinite if and only if $\alpha  \geq {\alpha _0} =  \frac{{{-\lambda _{\min }}(RD(G))}}{{k - {\lambda _{\min }}(RD(G))}}.$
\begin{proof}
Clearly, $R{D_\alpha }(G) = \alpha kI_{n} + (1 - \alpha )RD(G)$ because $G$ is $k$-reciprocal transmission regular. This implies that
${\alpha }k + (1 - {\alpha }){\lambda _{\min}}(RD(G)) = {\lambda _{\min}}(R{D_{{\alpha }}}(G))$.
So, $ {\lambda _{\min}}(R{D_{{\alpha }}}(G))\geq 0$ if and only if $\alpha\geq{\alpha _0} = \frac{{ - {\lambda _{min}}(RD(G))}}{{k - {\lambda _{min}}(RD(G))}}$, as required.
\end{proof}

\paragraph{Theorem 3.11.} For a complete bipartite graph $K_{a,n-a}$ of order $n\geq 4$ with $1\leq a\leq \lfloor\frac{n}{2}\rfloor$, the matrix $R{D_\alpha }(K_{a,n-a})$ is a  positive semidefinite matrix if and only if $\alpha\in[\alpha_{0}, 1]$, where the smallest $\alpha_{0}= \frac{n-1+3a(n-a)}{2n(n-1)+4a(n-a)}$.
\begin{proof}
In accordance with Corollary 2.3, the eigenvalues of $R{D_\alpha }(K_{a,n-a})$ are
\begin{itemize}
  \item $\frac{{ (2n - a)\alpha - 1}}{2}$, repeated $a-1$ times,
  \item $\frac{{(n + a)\alpha - 1}}{2}$, repeated $n-a-1$ times,
  \item the remaining two eigenvalues
  \[
  \frac{{(\alpha  + \frac{1}{2})n - 1 \pm \sqrt {{{(\alpha  - \frac{1}{2})}^2}{{(2a - n)}^2} + 4{{(1 - \alpha )}^2}a(n - a)} }}{2}.
  \]
\end{itemize}
Observe that $\frac{{(n + a)\alpha - 1}}{2}\leq\frac{{ (2n - a)\alpha - 1}}{2}$ as $a\leq \lfloor\frac{n}{2}\rfloor$. Thus, the minimum eigenvalue of $R{D_\alpha }(K_{a,n-a})$ is $\frac{{(n + a)\alpha - 1}}{2}$ or $\frac{{(\alpha  + \frac{1}{2})n - 1 -\sqrt {{{(\alpha  - \frac{1}{2})}^2}{{(2a - n)}^2} + 4{{(1 - \alpha )}^2}a(n - a)} }}{2}$. Consider the following functions
$\phi(x) = \frac{{(n + a)x - 1}}{2}$ and
\begin{equation*}
\varphi(x) = \frac{{(x + \frac{1}{2})n - 1 - \sqrt {{{(x  - \frac{1}{2})}^2}{{(2a - n)}^2} + 4{{(1 - x )}^2}a(n - a)} }}{2}
\end{equation*}
for $0\leq x\leq1$. By a simple calculation, the zero of the function $\varphi(x)$ is ${x_0} = \frac{n-1+3a(n-a)}{2n(n-1)+4a(n-a)}$, but $\phi(x_0)>0$. On the other hand, Proposition 3.5 implies that $\phi(x)$ and $\varphi(x)$ are strictly increasing in $x$. Therefore, ${x_0}$ is the smallest value ${\alpha _0}$ for $\alpha\in[0,1]$ such that $R{D_\alpha}(K_{a,n-a})$ is a positive semidefinite matrix.
\end{proof}

\paragraph{Theorem 3.12} Let $W(n)$ be a wheel graph, where $n\geq 4$. Then the matrix $R{D_\alpha }(W(n))$ is a positive semidefinite matrix if and only if $\alpha\in[\alpha_{0}, 1]$, where the smallest $\alpha_{0}=\frac{3}{{n + 5}}$ for $n=2k+1$, or $$\alpha_{0}=\frac{{1 - 2\cos (\frac{{2k\pi }}{{2k + 1}})}}{{n + 3 - 2\cos (\frac{{2k\pi }}{{2k + 1}})}},\;\;\;\; \text{for} \;\;n=2k+2.$$

\begin{proof}
From Corollary 2.5, the eigenvalues of the matrix $R{D_\alpha }(W(n))$ are
\[\frac{{\alpha (n + 3) - 1 + 2(1 - \alpha )\cos (\frac{{2j\pi }}{{n - 1}})}}{2}\;\; \text{for} \;\;1\leq j\leq n-2,
\]
and
\[
\frac{{(\alpha  + \frac{1}{2})n \pm \sqrt {{{(\alpha (n - 2) - \frac{n}{2})}^2} + 4{{(1 - \alpha )}^2}(n - 1)} }}{2}.
\]

Now we need to consider the minimum eigenvalue of the matrix $R{D_\alpha }(W(n))$. First we suppose that $n=2k+1$. Then the candidates to be the minimum eigenvalue of $R{D_\alpha }(W(n))$ are $\frac{\alpha (n + 5) -3}{2}$ and $\frac{{(\alpha  + \frac{1}{2})n - \sqrt {{{(\alpha (n - 2) - \frac{n}{2})}^2} + 4{{(1 - \alpha )}^2}(n - 1)} }}{2}$. Set $\phi(x) = \frac{(n + 5)x -3}{2}$ and
$$
\varphi(y) =\frac{{(y + \frac{1}{2})n - \sqrt {{{( (n - 2)y - \frac{n}{2})}^2} + 4{{(1 - y )}^2}(n - 1)} }}{2}.
$$
After a simple calculation, the zero of the function $ \phi(x) $ is ${x_0} = \frac{3}{{n + 5}}$, but $\varphi(x_0)>0$. Proposition 3.5 implies that $\phi(x)$ and $\varphi(x)$ are strictly increasing in $x$. Therefore, ${x_0}$ is the smallest value ${\alpha _0}$ for $\alpha\in[0,1]$, such that $R{D_\alpha }(W(n))$ is a positive semidefinite matrix.

Next assume that $n=2k+2$. We can find that the eigenvalue $\frac{{\alpha (n + 3) - 1 + 2(1 - \alpha )\cos (\frac{{2k\pi }}{{2k+ 1}})}}{2}$ or $\frac{{(\alpha  + \frac{1}{2})n - \sqrt {{{(\alpha (n - 2) - \frac{n}{2})}^2} + 4{{(1 - \alpha )}^2}(n - 1)} }}{2}$ will be the minimum eigenvalue of  $R{D_\alpha }(W(n))$. Set
$$
\phi(x) = \frac{{ (n + 3)x - 1 + 2(1 - x )\cos (\frac{{2k\pi }}{{2k+ 1}})}}{2}
$$
and
$$
\varphi(y) =\frac{{(y + \frac{1}{2})n - \sqrt {{{( (n - 2)y - \frac{n}{2})}^2} + 4{{(1 - y )}^2}(n - 1)} }}{2}.
$$
After a simple calculation, the zero of the function $ \phi(x) $ is ${x_0} = \frac{{1 - 2\cos (\frac{{2k\pi }}{{2k + 1}})}}{{n + 3 - 2\cos (\frac{{2k\pi }}{{2k + 1}})}}$, and the zero of the function $ \varphi(y)$ is ${y_0}= \frac{2}{{n + 4}}$. It is clear that $x_0\geq y_0$ whenever $n=2k+2$. From Proposition 3.5, $\phi(x)$ and $\varphi(x)$ are strictly increasing in $x$. Therefore, ${x_0}$ is the smallest value ${\alpha _0}$ for $\alpha\in[0,1]$, such that $R{D_\alpha }(W(n))$ is a positive semidefinite matrix.
\end{proof}

Next, we shall give some lower and upper bounds on the spectral radius of $R{D_\alpha }(G)$ in terms of the reciprocal transmission sequence $\{RTr_G(v_{1}),RTr_G(v_{2}),\ldots,RTr_G(v_{n})\}$ of a graph $G$.

\paragraph{Lemma 3.13.\cite{You}} Let $M=(m_{i,j})_{n\times n}$ be a complex irreducible matrix and $l_k=|\{m_{k,j}: m_{k,j}\neq 0, j\in \{1,2,\ldots,n\}\setminus\{k\}\}|$. Then
\begin{equation}\label{11}
\rho (M) \leq \mathop {\max }\limits_{1 \leq i \leq n} \left\{ {\left| {{m_{i,i}}} \right| + \sqrt {\sum\limits_{k = 1,k \ne i}^n {{l_k}{{\left| {{m_{k,i}}} \right|}^2}} } } \right\}.
\end{equation}
Moreover, if the equality holds in (\ref{11}), then for $\forall\;  i,j\in\{1,2,\ldots,n\}$,
\[
{\left| {{m_{i,i}}} \right| + \sqrt {\sum\limits_{k = 1,k \ne i}^n {{l_k}{{\left| {{m_{k,i}}} \right|}^2}} } }={\left| {{m_{j,j}}} \right| + \sqrt {\sum\limits_{k = 1,k \ne j}^n {{l_k}{{\left| {{m_{k,j}}} \right|}^2}} } }.
\]

Now, applying Lemma 3.13 to $RD_{\alpha}(G)$, we may conclude the following.

\paragraph{Theorem 3.14.} Let $G$ be a connected graph with $n$ vertices. Then
\begin{equation}\label{12}
\rho (R{D_\alpha }(G)) \leq \mathop {\max }\limits_{1 \leq i \leq n} \left\{ {\alpha RT{r_G}(v_{i}) + (1-\alpha)\sqrt {(n - 1)\sum\limits_{k = 1,k \ne i}^n {{{\left( {\frac{1}{{{d_{G}(v_{k},v_{i})}}}} \right)}^2}} } } \right\}.
\end{equation}
If the equality holds in (\ref{12}), then
$\alpha RT{r_G}({v_i}) + (1 - \alpha )\sqrt {(n - 1)\sum\limits_{k = 1,k \ne i}^n {{{\left( {\frac{1}{{{d_G}({v_k},{v_i})}}} \right)}^2}} }$
are all equal for any $ i\in\{1,2,\ldots,n\}$.

\paragraph{Lemma 3.15.\cite{Minc}} Let $M=(m_{i,j})_{n\times n}$ be a complex irreducible matrix. Then
\begin{equation}\label{13}
 \mathop {\min }\limits_{1 \leq i \leq n} \left\{ {\sum\limits_{j = 1}^n {{m_{i,j}}} } \right\} \leq \rho (M) \leq \mathop {\max }\limits_{1 \leq i \leq n} \left\{ {\sum\limits_{j = 1}^n {{m_{i,j}}} } \right\}.
\end{equation}
Moreover, each one of the equalities holds in (\ref{13}) if and only if all row sums of $M$ are equal.

\paragraph{Theorem 3.16.} Let $G$ be a connected graph of order $n$. Then
\begin{equation}\label{14}
\begin{gathered}
  \mathop {\min }\limits_{1 \leq i \leq n} \left\{ {\alpha RT{r_G}\left( v_{i} \right) + \frac{{(1 - \alpha )R{T_i}}}{{RT{r_G}(v_{i})}}} \right\} \leq \rho (R{D_\alpha }(G)) \leq \mathop {\max }\limits_{1 \leq i \leq n} \left\{ {\alpha RT{r_G}\left( v_{i} \right)+\frac{{(1 - \alpha )R{T_i}}}{{RT{r_G}(v_{i})}}} \right\},
\end{gathered}
\end{equation}
where $R{T_i} = \sum\nolimits_{j = 1,j \ne i}^n {R{D_{i,j}}RT{r_G}({v_j})} $.
Moreover, for $\frac{1}{2} \leq \alpha  \leq 1$, the equalities hold in (\ref{14}) if and only if graph $G$ is reciprocal transmission regular.
\begin{proof}
Let ${(RT(G))^{-1}}$ be the inverse matrix of $RT(G)$. By a simple calculation, the $i$-th row sum of ${(RT(G))^{ - 1}}R{D_\alpha }(G)RT(G)$ is equal to
\[
\alpha {RT{r_G}(v_{i})} + (1 - \alpha )\frac{{R{T_i}}}{{RT{r_G}(v_{i})}}.
\]
Use Lemma 3.15 by taking $M ={(RT(G))^{ - 1}}R{D_\alpha }(G)RT(G)$, the inequality (\ref{14}) is immediate.

Assume that graph $G$ is  $k$-reciprocal transmission regular. It is easy to verify that $\rho (R{D_\alpha }(G)) = k$ and $\alpha RTr_G(v) + (1 - \alpha )\frac{{R{T_v}}}{{RT{r_G}(v)}} = k$ for any $v\in V(G)$. Hence, the equality holds in (\ref{14}).

Conversely, assume either of the equalities holds in (\ref{14}). At that time, Lemma 3.15 implies that all row sums of ${(RT(G))^{ - 1}}R{D_\alpha }(G)RT(G)$ are equal. Now, for $\forall \;u,v\in V(G)$, one obtains that
\begin{equation}\label{15}
{\alpha RT{r_G}\left( u \right) + (1 - \alpha )\frac{{R{T_u}}}{{RT{r_G}(u)}}}={\alpha RT{r_G}\left( v \right) + (1 - \alpha )\frac{{R{T_v}}}{{RT{r_G}(v)}}}.
\end{equation}
Without loss of generality, let $RT{r_G}\left( u \right) = RT{r_{\max }}$ and $RT{r_G}\left( v \right) = RT{r_{\min }}$, then we have $RT\left( u \right) \geq RT{r_{\max }}RT{r_{\min }}$, $RT\left( v \right) \leq RT{r_{\max }}RT{r_{\min }}$. According to (\ref{15}), it is easy to see that
\[
\alpha RT{r_{\max }} + (1 - \alpha )RT{r_{\min }} \leq \alpha RT{r_{\min }} + (1 - \alpha )RT{r_{\max }}.
\]
This implies that $RT{r_{\min}} = RT{r_{\max}}$ for $\frac{1}{2}\leq\alpha<1$. Hence, the graph $G$ is reciprocal transmission regular.
\end{proof}

\paragraph{Lemma 3.17.\cite{Zhou}} Let $M=(m_{i,j})_{n\times n}$ be a nonnegative irreducible symmetric matrix  with the $i$-th row sum $M_{i}$. Then
\begin{equation}\label{16}
 \sqrt {\frac{{\sum\limits_{j = 1}^n {M_{i}^2} }}{n}}  \leq \rho (M) \leq \mathop {\max }\limits_{1 \leq i \leq n} \sum\limits_{j = 1}^n {{m_{i,j}}} \sqrt {\frac{{{M_j}}}{{{M_i}}}}.
\end{equation}
Each one of the equalities in (\ref{16}) holds if and only if either all row sums of $M$ are equal, or there exists a permutation matrix $Q$ such that
\[
{Q^{\rm T}}MQ = \left( {\begin{array}{*{20}{c}}
  0&C \\
  {{C^{\rm T}}}&0
\end{array}} \right),
\]
where all row sum of $C$ are equal.

\paragraph{Theorem 3.18.}
Let $G$ be a connected graph of order $n$, where $n\geq2$. Then
\begin{equation}\label{17}
\sqrt {\frac{{\sum\limits_{i = 1}^n {RT{r_G}{{\left( v_{i} \right)}^2}} }}{n}}  \leq \rho (R{D_\alpha }(G)) \leq \mathop {\max }\limits_{1 \leq i \leq n} (\alpha RT{r_G}(v_{i}) + (1 - \alpha )\sum\limits_{j = 1,j \ne i}^n {\frac{1}{{{d_{G}(v_i,v_j)}}}} \sqrt {\frac{{RT{r_G}(v_j)}}{{RT{r_G}(v_i)}}} ).
\end{equation}
Either of the equalities holds in (\ref{17}) if and only if graph $G$ is reciprocal transmission regular.
\begin{proof}
Since $R{D_\alpha }(G)$ is a nonnegative irreducible symmetric matrix. Then, it follows from Lemma 3.17 that the inequality (\ref{17}) is immediate. Clearly, there is not a permutation matrix $Q$ such that \[
{Q^{\rm T}}R{D_\alpha }(G)Q = \left( {\begin{array}{*{20}{c}}
	0&C \\
	{{C^{\rm T}}}&0
	\end{array}} \right).
\]
So, the required result follows.
\end{proof}

Applying Theorem 3.18, we easily obtain the following result.

\paragraph{Corollary 3.19.}
Let $G$ be a connected graph with $n\geq2$ vertices and Harary index $H(G)$. Then
\begin{equation}\label{18}
\rho (R{D_\alpha }(G))\geq\frac{{2H(G)}}{n}.
\end{equation}
The equality holds in (\ref{18}) if and only if graph $G$ is  reciprocal transmission regular.
\begin{proof}
Using the Cauchy-Schwarz inequality to the left-side of the inequality (\ref{17}), we have
\[\rho (R{D_\alpha }(G)) \geq \sqrt {\frac{{\sum\limits_{i = 1}^n {RT{r_G}{{\left( v_i \right)}^2}} }}{n}}  \geq \frac{{\sum\limits_{i = 1}^n {RT{r_G}\left( v_i \right)} }}{n} = \frac{{2H(G)}}{n}.
\]
Clearly, the equality holds if and only if $G$ is reciprocal transmission regular.
\end{proof}

Remark that, if $\alpha=0$, then this result reduces to $\rho (R{D}(G))\geq\frac{{2H(G)}}{n}$. Similarly, let $\alpha=\frac{1}{2}$, then this result reduces to $\rho (RQ(G))\geq\frac{{4H(G)}}{n}$. So, Corollary 3.19 generalizes the results in \cite{Alhevaz} and \cite{Zhou}.

\section{Graphs with fixed structure parameters}

Given a connected graph $G=(V(G),E(G))$, its vertex connectivity $\kappa(G)$ is the minimum number of vertices whose removal yields a disconnected graph. An edge cut $[S, \overline{S}]$ of $G$ is a subset of $E(G)$ between $S$ and $\overline{S}$, where $S$ is a nonempty proper subset of $V(G)$ and $\overline{S}=V(G)\backslash {S}$. A $k$-edge cut is an edge cut of $k$ edges in $G$. The edge connectivity $\kappa'(G)$ of $G$ is the minimum cardinal of the edge cuts in $G$. Obviously, for the minimum vertex degree $\delta (G)$, vertex connectivity $\kappa(G)$ and edge connectivity $\kappa'(G)$ in the graph $G$, we know that $\kappa (G) \leq \kappa' (G) \leq \delta (G)$.

Let $G_{n}^{r}$ and $\overline{G}_{n}^{r}$ denote the set of all graphs of order $n$ with vertex connectivity $r$ and edge connectivity $r$, respectively. Obviously, $G_n^{n - 1} = \overline G_n^{n - 1} = {K_n}$.

Let $R{D_\alpha }(G)$ be the $R{D_\alpha }$-matrix of graph $G$ and $\textbf{x}=(x_1,x_2,\ldots,x_n)^{\rm T}$ be a column vector of order $n$. As we all know,
\begin{equation*}
{\textbf{x}^{\rm T}}R{D_\alpha }(G)\textbf{x} = \alpha \sum\limits_{i = 1}^n {RT{r_G}({v_i})} {x_i}^2 + 2(1 - \alpha )\sum\limits_{1 \leq i < j \leq n}^{} {\frac{1}{{{d_G}({v_i},{v_j})}}} {x_i}{x_j}.
\end{equation*}
Then
\begin{equation}\label{19}
  \begin {aligned}
  {\textbf{x}^{\rm T}}R{D_\alpha }(G')\textbf{x} - {\textbf{x}^{\rm T}}R{D_\alpha }(G)\textbf{x} = &\alpha\sum\limits_{i = 1}^n {(RT{r_{G'}}({v_i}) - RT{r_G}({v_i}))} {x_i}^2 \\
  &+ 2(1 - \alpha )\sum\limits_{1 \leq i < j \leq n}^{} {(\frac{1}{{{d_{G'}}({v_i},{v_j})}} - \frac{1}{{{d_G}({v_i},{v_j})}}} ){x_i}{x_j}.
  \end {aligned}
\end{equation}

\paragraph{Theorem 4.1.} Let $n$ and $r$ be given  integers with $1\leq r\leq n-2$. Then, for $0\leq\alpha<1 $, the graph ${K_r} \vee ({K_1} \cup {K_{n - r - 1}})$ is the unique graph with maximal spectral radius of the matrix $RD_{\alpha}(G)$ in $G\in G_{n}^{r}$.
\begin{proof}
Assume that $G$ is a graph with maximal spectral radius $\rho (R{D_\alpha }(G))$ of $RD_{\alpha}(G)$ in $G\in G_{n}^{r}$. Proposition 3.2 implies that the graph $G$ must be isomorphic to ${K_r} \vee ({K_{{n_1}}} \cup {K_{{n_2}}})$, where $n_{1}+n_{2}=n-r$. Without loss of generality, suppose that $1\leq n_{1}\leq n_{2}$.

Suppose towards contradiction that $n_{1}>1$. Let $\textbf{x}$ be the perron eigenvector with respect to the spectral radius $\rho (R{D_\alpha }(G))$ of $RD_{\alpha}(G)$ and $\left\| \textbf{x} \right\|=1$. According to the equivalence of vertices in $G$, the perron eigenvector $\textbf{x}$ can be written as
\[{{\mathbf{x}}^{\text{T}}} = (\underbrace {{x_1}, \ldots ,{x_1}}_{{n_1}},\underbrace {{x_2}, \ldots ,{x_2}}_{{n_2}},\underbrace {{x_3}, \ldots ,{x_3}}_r).\]
Let $v_{1}\in V(K_{n_1})$. Construct a new graph $G'$  from graph $G$ by deleting the edges $\{{v_1}{v_i}:{v_i} \in V({K_{n_1}})\backslash \{v_1\}\}$ and adding the edges $\{{v_1}{v_j}:{v_j} \in V({K_{n_2}})\}$, that is,
\[G' = G -\bigcup\limits_{{v_i} \in V({K_{n_1}})\backslash \{v_1\}} v_{1}v_{i} + \bigcup\limits_{{v_j} \in V({K_{n_2}})} {{v_1}{v_j}}.
\]
Obviously, $G'= {K_r} \vee ({K_{{n_1-1}}} \cup {K_{{n_2}+1}})$. Since
\begin{equation*}
  \begin {aligned}
  \rho (R{D_\alpha }(G') - \rho (R{D_\alpha }(G)) \geq {\textbf{x}^{\rm T}}R{D_\alpha }(G')\textbf{x} - {\textbf{x}^{\rm T}}R{D_\alpha }(G)\textbf{x}.
  \end {aligned}
\end{equation*}
According to (\ref{19}), we only need to compute ${\frac{1}{{{d_{G'}}({v_i},{v_j})}} - \frac{1}{{{d_G}({v_i},{v_j})}}}$ and ${RT{r_{G'}}({v_i}) - RT{r_G}({v_i})}$ for any  $v_{i},v_{j}\in V(G)$. After some careful observations and calculations, one gets that

\begin{enumerate}[(i)]
  \item ${\frac{1}{{{d_{G'}}({v_1},{v_i})}} - \frac{1}{{{d_G}({v_1},{v_i})}}}=-\frac{1}{2}$
  for any $ {v_i}\in V(K_{n_{1}})\backslash \{v_{1}\}$;
  \item ${\frac{1}{{{d_{G'}}({v_1},{v_j})}} - \frac{1}{{{d_G}({v_1},{v_j})}}}=\frac{1}{2}$
  for any ${v_j}\in V(K_{n_{2}})$;
  \item $RT{r_{G'}}({v_1}) - RT{r_G}({v_1}) = \frac{{{n_2} - ({n_1} - 1)}}{2}$;
  \item $RT{r_{G'}}({v_i}) - RT{r_G}({v_i}) =  - \frac{1}{2}$ for any  $ {v_i}\in V(K_{n_{1}})\backslash\{v_{1}\}$;
  \item $RT{r_{G'}}({v_j}) - RT{r_G}({v_j}) = \frac{1}{2}$
  for any  ${v_j}\in V(K_{n_{2}})$.
\end{enumerate}
Plugging the above into (\ref{19}), we have
\[
\begin {aligned}
{\textbf{x}^{\rm T}}R{D_\alpha }(G')\textbf{x} - {\textbf{x}^{\rm T}}R{D_\alpha }(G)\textbf{x}=&\alpha \left[\frac{{{n_2} - ({n_1} - 1)}}{2}{x_1^2} - \frac{{({n_1} - 1)}}{2}{x_1^2} + \frac{{{n_2}}}{2}{x_2^2}\right]\\
&+ (1 - \alpha )[{n_2}{x_1}{x_2} - ({n_1} - 1){x_1^2}]\\
= &\frac{{\alpha {n_2}}}{2}({x_1^2} + {x_2^2}) - ({n_1} - 1){x_1^2} + (1 - \alpha ){n_2}{x_1}{x_2}.
\end {aligned}
\]
As a matter of convenience, let $\rho=\rho (R{D_\alpha }(G))$. It follows from $R{D_\alpha }(G)\textbf{x} = \rho \textbf{x}$ that
\[
\rho {x_1} = \alpha (n - 1 - \frac{{{n_2}}}{2}){x_1} + (1 - \alpha )[({n_1} - 1){x_1} + \frac{{{n_2}}}{2}{x_2} + r{x_3}]
\]
and
\[
\rho {x_2} = \alpha (n - 1 - \frac{{{n_1}}}{2}){x_2} + (1 - \alpha )[\frac{{{n_1}}}{2}{x_1} + ({n_2} - 1){x_2} + r{x_3}],
\]
or equivalently,
\[
\rho {x_1} - \alpha (n - 1 - \frac{{{n_2}}}{2}){x_1} - (1 - \alpha )({n_1} - 1){x_1} - (1 - \alpha )\frac{{{n_2}}}{2}{x_2} = (1 - \alpha )r{x_3}
\]
and
\[
\rho {x_2} - \alpha (n - 1 - \frac{{{n_1}}}{2}){x_2} - (1 - \alpha )({n_2} - 1){x_2} - (1 - \alpha )\frac{{{n_1}}}{2}{x_1} = (1 - \alpha )r{x_3}.
\]
By combining the above two equations, we easily obtain that
\begin{equation}\label{20}
(\rho  - \alpha n + \frac{{{n_1} + {n_2}}}{2}\alpha  + 1 - \frac{{{n_1}}}{2}){x_1} = (\rho  - \alpha n + \frac{{{n_1} + {n_2}}}{2}\alpha  + 1 - \frac{{{n_2}}}{2}){x_2}.
\end{equation}
Since $(1 - \alpha )\frac{{{n_1}}}{2}{x_1} + (1 - \alpha )\frac{{{n_2}}}{2}{x_2} > 0$, then $(\rho  - \alpha n + \frac{{{n_1} + {n_2}}}{2}\alpha  + 1 - \frac{{{n_1}}}{2}){x_1}>0$.
It follows from (\ref{20}) that $\frac{{{x_2}}}{{{x_1}}} \geq1$ as $n_{2}\geq n_{1}$. Thus,
\[
\frac{{\alpha {n_2}}}{2}({x_1^2} + {x_2^2}) - ({n_1} - 1){x_1^2} + (1 - \alpha ){n_2}{x_1}{x_2}={x_1^2}[\frac{{{\alpha n_2}}}{2}{(\frac{{{x_2}}}{{{x_1}}} - 1)^2} + ({n_2}\frac{{{x_2}}}{{{x_1}}} - {n_1} + 1)] > 0,
\]
which leads to $\rho (R{D_\alpha }(G')) > \rho (R{D_\alpha }(G))$. This contradicts our previous assumption. Hence, $G={K_r} \vee ({K_1} \cup {K_{n - r - 1}})$.
\end{proof}

In what follows, we intends to give the extremal graphs with maximal spectral radius of the $RD_{\alpha}$-matrix among all connected graphs of fixed edge connectivity. To do this, we first recall the following lemma.

\paragraph{Lemma 4.2.\cite{West}}
Let $S$ be a nonempty vertex subset of a graph $G$ with minimum vertex degree $\delta (G)$. If $\left| {\left[ {S,\overline S } \right]} \right| < \delta (G)$, then $\left| S \right| > \delta (G)$.

\paragraph{Theorem 4.3.} For $1\leq r\leq n-2$ and $0\leq\alpha<1 $, the graph ${K_r} \vee ({K_1} \cup {K_{n - r - 1}})$ is the unique graph with maximal spectral radius of the matrix $RD_{\alpha}(G)$ in $ G\in \overline G_{n}^{r}$.
\begin{proof}
Let $G$ be the graph with maximal spectral radius of the matrix $RD_{\alpha}(G)$ in $ G\in \overline G_{n}^{r}$. It is known that $\delta (G) \geq r$. If there exists a vertex $v$ in graph $G$ with degree $r$, then $[ \{v\},V(G)\backslash \{v\}]$ is an $r$-edge cut of $G$. Proposition 3.2 implies that the induced subgraph of $V(G)\backslash \{v\}$ must be a complete graph. So, the required result follows.

Now assume that $\delta(G)>r$ and ${\left[ {S,\overline S } \right]}$ is an $r$-edge cut of $G$ with $|S|=n_1$ and $|\overline S|=n_2$. Let ${G_1} = G[S]$ and ${G_2} = G[\overline S]$ be two induced subgraphs of $G$ with respect to $S$ and $\overline{S}$, respectively.  It follows from Proposition 3.2 that $G_{1}$ and $G_{2}$ are complete graphs. Noting that $\delta(G)>r$, then ${n_1} > 1$ and ${n_2} > 1$. Let $V({G_1}) = \{{{v_1},{v_2},\ldots,{v_{{n_1}}}}\}$ and $V({G_2}) = \{v_{{n_1} + 1},{v_{{n_1} + 2}},\ldots,{v_{n}}\}$. Also let $ \textbf{x}= ({x_1}, {x_2}\ldots,x_{n_{1}}, x_{n_{1}+1},x_{n_{1}+2}, \ldots ,x_{n})$ be the perron eigenvector with respect to $\rho(RD_{\alpha}(G))$, where ${x_i}$ be the component of $\textbf{x}$ corresponding to the vertex $v_{i}$ in $V(G)$. Without loss of generality, let ${x_1}=\mathop {\min}\nolimits_{1\leq i\leq n} \{{x_i}\}$ and $x_{1}\leq x_{2}\leq\cdots \leq x_{n_{1}}$. We may assume that $v_{1}$ is adjacent to $t$ vertices of $G_{2}$. Obviously, $t\leq \min\{r, n_{2}\}$.

Now consider the following two cases: $t=r$ and $t<r$.

\textit{Case 1:} $t=r$. In this case, the edges between $v_{1}$ and $G_{2}$ are exactly the cut edges between $G_{1}$ between $G_{2}$. At that time, we have ${n_2} \geq r + 2$. Indeed, if $n_{2}=r$, or $n_{2}=r+1$, then there must exists a vertex with degree $r$ in $G_{2}$. This contradicts our hypothesis $\delta(G)>r$. Now construct a new graph $G'$ from $G$ as follows:
\[
G' = G - \bigcup\limits_{{v_i} \in V({G_1})\setminus{v_{1}}} {{v_1}{v_i}}  + \bigcup\limits_{{v_i} \in V({G_1})\backslash \{v_1\},\;{v_j} \in V({G_2})} {{v_i}{v_j}}.
\]
It is clear that $G' = {K_r} \vee ({K_{1}} \cup {K_{n-r-1}})$.
Let $A=V(G_{1})\backslash\{v_{1}\}$, $B=\{v_j\in V(G_2): v_1v_j\in E(G)\}$ and $C=V(G_{2})\backslash B$.
According to (\ref{19}), we only need to calculate
${\frac{1}{{{d_{G'}}({v_i},{v_j})}} - \frac{1}{{{d_G}({v_i},{v_j})}}}$ and ${RT{r_{G'}}({v_i}) - RT{r_G}({v_i})}$ for any $v_{i},v_{j}\in V(G)$. By some observations and calculations, one has
\begin{enumerate}[(i)]
  \item ${\frac{1}{{{d_{G'}}({v_1},{v_i})}} - \frac{1}{{{d_G}({v_1},{v_i})}}}=-\frac{1}{2}$
  for any $ {v_i}\in A$;
  \item ${\frac{1}{{{d_{G'}}({v_i},{v_j})}} - \frac{1}{{{d_G}({v_i},{v_j})}}}=\frac{1}{2}$
  for any ${v_i}\in A$, ${v_j}\in B$;
  \item ${\frac{1}{{{d_{G'}}({v_i},{v_k})}} - \frac{1}{{{d_G}({v_i},{v_k})}}}=\frac{2}{3}$
  for any ${v_i}\in A$, $ {v_k}\in C$;
  \item $RT{r_{G'}}({v_1}) - RT{r_G}({v_1}) = \frac{{ - |A|}}{2}$;
  \item $RT{r_{G'}}({v_i}) - RT{r_G}({v_i}) = - \frac{1}{2} + \frac{1}{2}|B| + \frac{2}{3}|C|$ for any $ {v_i}\in A$;
  \item $RT{r_{G'}}({v_j}) - RT{r_G}({v_j}) = \frac{1}{2}|A|$
  for any ${v_j}\in B$;
 \item $RT{r_{G'}}({v_k}) - RT{r_G}({v_k}) = \frac{2}{3}|A|$
  for any $ {v_k}\in C$.
\end{enumerate}
Now let us substitute the above into (\ref{19}),
\[\small
\begin {aligned}
\rho (R{D_\alpha }(G') - \rho (R{D_\alpha }(G)) \geq& {\textbf{x}^{\rm T }}R{D_\alpha }(G')\textbf{x} - {\textbf{x}^{\rm T }}R{D_\alpha }(G)\textbf{x}\\
 =& \alpha \left [ - \frac{1}{2}\left| A \right|{x_1^2} + ( - \frac{1}{2} + \frac{1}{2}\left| B \right| + \frac{2}{3}\left| C \right|)\sum\limits_{{v_i} \in A} {{x_i^2}}  + \frac{1}{2}\left| A \right|\sum\limits_{{v_j} \in B} {{x_j^2}}+ \frac{2}{3}\left| A \right|\sum\limits_{{v_k} \in C} {{x_k^2}} \right ] \\
&+ (1 - \alpha )\left[ - \sum\limits_{{v_i} \in A} {{x_1}} {x_i} + \sum\limits_{{v_i} \in A,{v_j} \in B} {{x_i}} {x_j} + \sum\limits_{{v_i} \in A,{v_k} \in C} {\frac{4}{3}{x_i}} {x_k}\right]\\
>&0,
\end {aligned}
\]
where the last inequality holds because ${x_1}=\mathop {\min}\nolimits_{1\leq i\leq n} \{{x_i}\}$. So, $\rho (R{D_\alpha }(G')) > \rho (R{D_\alpha }(G))$, which yields a contradiction.

\textit{Case 2:} $t<r$. In this case, Lemma 4.2 implies that $n_{2}>r+1$ as $\delta (G)>r$. Let $W = \left\{ {v \in V({G_2}):vw\notin E(G), w \in V({G_1})\backslash \{v_1\}} \right\}$. Since $|[V({G_1})\backslash \{v_1\},V(G_{2})]|=r-t$ and $n_{2}>r+1$, then $| W| >0$. Let ${V_{F}} = \{{v_2}, \ldots ,{v_{{n_1}-r+ t}}\}\subseteq V(G_{1})\backslash \{v_{1}\}$. For $\forall v_{f}\in V_{F}$, $\forall v_{j}\in G_{2}$, vertex $v_{f}$ is not adjacent to $v_{j}$. Let ${V_{F'}} =V(G_{1})\setminus \{\{v_{1}\}\cup {V_{F}}\}$. It is clear that $|[{V_{F'}}, V(G_2)]|=r-t$. Now construct a new graph $G''$ from $G$ as below:
\[
G''= G - \bigcup\limits_{{v_f} \in V(F)} {{v_1}{v_f}}  + \bigcup\limits_{{v_i} \in V({G_1})\backslash \{v_1\},\;{v_j} \in V({G_2}),\; v_iv_j\notin E(G)} {{v_i}{v_j}}.
\]
Obviously, $G'' = {K_r} \vee ({K_{1}} \cup {K_{n-r-1}})$. In the following, we calculate
${\frac{1}{{{d_{G''}}({v_i},{v_j})}} - \frac{1}{{{d_G}({v_i},{v_j})}}}$ and ${RT{r_{G''}}({v_i}) - RT{r_G}({v_i})}$ for any $v_{i},v_{j}\in V(G)$.
\begin{enumerate}[(i)]
  \item ${\frac{1}{{{d_{G''}}({v_1},{v_f})}} - \frac{1}{{{d_G}({v_1},{v_f})}}}=-\frac{1}{2}$
  for any $ {v_f}\in {V_{F}}$;
  \item ${\frac{1}{{{d_{G''}}({v_i},{v_j})}} - \frac{1}{{{d_G}({v_i},{v_j})}}}=1-\frac{1}{d_{G}(v_i,v_j)}$
  for any $ {v_i}\in V(G_{1})\setminus \{v_{1}\}$, ${v_j}\in V(G_{2})$;
  \item $RT{r_{G''}}({v_1}) - RT{r_G}({v_1}) = \frac{{ - |F|}}{2}$;
  \item $a_{1}=RT{r_{G''}}({v_f}) - RT{r_G}({v_f}) = - \frac{1}{2} + \sum\limits_{{v_j} \in V({G_2})} {(1 - \frac{1}{{{d_G}(v_f,v_j)}})} $ for any $ {v_f}\in {V_{F}}$;
  \item $a_{2}=RT{r_{G''}}({v_{f'}}) - RT{r_G}({v_{f'}}) =\sum\limits_{ {v_j} \in V({G_2})} {(1 - \frac{1}{{{d_G}(v_{f'},v_j)}})} $
  for any $ {v_{f'}}\in {V_{F'}}$;
 \item $a_{3}=RT{r_{G''}}({v_j}) - RT{r_G}({v_j}) =\sum\limits_{ {v_i} \in V({G_1})\setminus \{v_1\}} {(1 - \frac{1}{{{d_G}(v_i,v_j)}})} $
  for any $  {v_j}\in V(G_{2}) $.
\end{enumerate}
Again, substitute the above into (\ref{19}),
\begin{equation}\label{21}
\small
  \begin {aligned}
  \rho (R{D_\alpha }(G'') - \rho (R{D_\alpha }(G)) &\geq {\textbf{x}^{\rm T} }R{D_\alpha }(G'')\textbf{x} - {\textbf{x}^{\rm T} }R{D_\alpha }(G)\textbf{x}\\
   &=\alpha\left[-\frac{|F|}{2}{x_1^2} + \sum\limits_{{v_f} \in {V_{F}}} {{a_1}{x_f^2}}  + \sum\limits_{{v_{f'}} \in {V_{F'}}} {{a_2}{x_{f'}^2}}  + \sum\limits_{{v_j} \in V({G_2})} {{a_3}{x_j^2}}\right] \\
&\;\;\;\;+(1 - \alpha )\left[ - \sum\limits_{{v_f} \in {V_{F}}} {{x_1}} {x_f} + \sum\limits_{{v_i} \in {G_1}\backslash \{v_1\},{v_j} \in V({G_2})} {2(1 - \frac{1}{{{d_G}(v_i,v_j)}}){x_i}} {x_j}\right].\\
  \end {aligned}
\end{equation}
Since every vertex in ${V_{F}}$ is not adjacent to any vertex of $G_{2}$, then ${{d_G}(v_f,v_j)}\geq 2$. This implies that $a_{1}>\frac{1}{2}(n_{2}-1)$.

Recall that $n_{2}>r+1$ and ${d_G}(v_f,v_j) \geq 2$ for any ${v_f} \in {V_F},{v_j} \in V({G_2})$. Thus,
\[
-\frac{|F|}{2}{x_1^2} + \sum\limits_{{v_f} \in {V_{F}}} {{a_1}{x_f^2}}>\frac{1}{2}|F|(n_{2}-2)x_{1}^{2}>0,
\]
and
\[\begin{gathered}\small
   {\sum _{{v_i} \in {G_1}\backslash \{v_1\},\;{v_j} \in V({G_2})}}2(1 - \frac{1}{{{d_G}(v_i,v_j)}}){x_i}{x_j}- \sum\limits_{{v_f} \in {V_F}} {{x_1}} {x_f}  \hfill \\
  \;\;\;\;\;=  \sum\limits_{{v_f} \in {V_F},{v_j} \in V({G_2})} {2(1 - \frac{1}{{{d_G}(v_f,v_j)}}){x_f}} {x_j} + \sum\limits_{{v_{f'}} \in {V_{F'}},{v_j} \in V({G_2})} {2(1 - \frac{1}{{{d_G}(v_{f'},v_j)}}){x_{f'}}} {x_j} - \sum\limits_{{v_f} \in {V_F}} {{x_1}} {x_f}  \hfill \\
    \;\;\;\;\;\geq  \sum\limits_{{v_f} \in {V_F},{v_j} \in V({G_2})} {{x_f}} {x_j} + \sum\limits_{{v_{f'}} \in {V_{F'}},{v_j} \in V({G_2})} {2(1 - \frac{1}{{{d_G}(v_{f'},v_j)}}){x_{f'}}} {x_j}- \sum\limits_{{v_f} \in {V_F}} {{x_1}} {x_f} \hfill \\
   \;\;\;\;\; > 0. \hfill \\
\end{gathered}
\]
It follows from (21) that $\rho (R{D_\alpha }(G'') - \rho (R{D_\alpha }(G))>0$, that is, $\rho (R{D_\alpha }(G'')) > \rho (R{D_\alpha }(G))$, a contradiction. To sum up, ${K_r} \vee ({K_1} \cup {K_{n - r - 1}})$  is the unique graph with maximal spectral radius of the matrix $RD_{\alpha}(G)$ in $ G\in \overline G_{n}^{r}$..
\end{proof}

For a connected graph $G$, the chromatic number $\chi(G)$ is the smallest number of colors of $V(G)$ such that any two adjacent vertices with different colors. The Tur\'{a}n graph, denoted by $T_{n,r}$, is the complete $r$-partite graph with $n$ vertices and each part has $\lfloor \frac{n}{r}\rfloor$ or $\lceil\frac{n}{r}\rceil$ vertices.
Denote the set of all connected graphs of order $n$ with chromatic number $\chi$ by $\widetilde{G}_{n}^{\chi}$.

\paragraph{Theorem 4.4.} For $n\geq 2$ and $0\leq \alpha \leq \frac{7}{16}$,
the Tur\'{a}n graph $T_{n,\chi}$ is the unique graph with maximal spectral radius of the matrix $RD_{\alpha}(G)$ in $G\in \widetilde{G}_{n}^{\chi}$.

\begin{proof}
Let $G$ be a graph with the maximal spectral radius of the matrix $RD_{\alpha}(G)$ in $G\in \widetilde{G}_{n}^{\chi}$. The conclusion is trivial when $\chi=n$, that is, $G$ is a complete graph $K_n$. Also noting that the unique connected non-complete graph is the path $P_3$ whenever $n=2$ or $3$. In the following, assume that $2\leq\chi\leq n-1$, $n\geq4$ and $G$ has a partition $\pi: V_{1}\cup V_{2}\cup\cdots\cup V_{\chi}$ of $V(G)$, where $|V_{i}|=n_{i}$ and $\sum\nolimits_{i = 1}^\chi  {{n_i}}=n $. Proposition 3.2 implies that $G=K_{n_{1},n_{2},\ldots,n_{\chi}}$. We may assume that $n_1=\max\nolimits_{1\leq i \leq \chi}\{n_{i}\}$ and $n_2=\min\nolimits_{1\leq i \leq \chi}\{n_{i}\}$. Next, we only need to prove that $n_1-n_2\leq1$.

Suppose towards contradiction that $n_{1}-n_{2}>1$. Let $\textbf{x}$ be the perron eigenvector corresponding to the spectral radius $\rho(RD_{\alpha}(G))$ of $RD_{\alpha}(G)$. According to the equivalence of the vertices in $G$, the vector $\textbf{x}$ can be written as
\[
\textbf{x}^{\rm T} = ({x_1}, \ldots ,{x_1},{x_2}, \ldots ,{x_2},\ldots,{x_\chi}, \ldots ,{x_\chi}).
\]
Now, select any vertex $v_{1}\in V_{1}$, construct graph $G'$ by deleting the edges $v_1v_j$ for any $v_{j}\in V_{2}$, and adding the adges $v_1v_i$ for any $v_{i}\in V_{1}\setminus \{v_1\}$ in graph $G$, that is,
\[
G' = G - \bigcup\limits_{{v_j} \in V_{2}} {{v_1}{v_j}}+\bigcup\limits_{{v_i} \in V_{1}\backslash {v_1}} {{v_1}{v_i}}.
\]
Obviously, $G' =K_{n_{1}-1,n_{2}+1,n_{3},\ldots,n_{\chi}}$.

In light of (\ref{19}), we first calculate
${\frac{1}{{{d_{G'}}({v_i},{v_j})}} - \frac{1}{{{d_G}({v_i},{v_j})}}}$ and ${RT{r_{G'}}({v_i}) - RT{r_G}({v_i})}$ for $v_{i},v_{j}\in V(G)$. These values are listed below:
\begin{enumerate}[(i)]
  \item ${\frac{1}{{{d_{G'}}({v_1},{v_i})}} - \frac{1}{{{d_G}({v_1},{v_i})}}}=\frac{1}{2}$
  for every ${v_i}\in V_{1}\backslash \{v_{1}\}$;
  \item ${\frac{1}{{{d_{G'}}({v_1},{v_j})}} - \frac{1}{{{d_G}({v_1},{v_j})}}}=-\frac{1}{2}$
  for every $ {v_j}\in V_{2}$;
  \item $RT{r_{G'}}({v_1}) - RT{r_G}({v_1}) = \frac{{{n_1}-{n_2} - 1}}{2}$;
  \item $RT{r_{G'}}({v_i}) - RT{r_G}({v_i}) = \frac{1}{2}$ for every $ {v_i}\in V_{1}\backslash\{v_{1}\}$;

  \item $RT{r_{G'}}({v_j}) - RT{r_G}({v_j})=-\frac{1}{2}$
  for every $ {v_j}\in V_{2}$.
\end{enumerate}
Substitute the above into (\ref{19}) to obtain
\begin{equation}\label{22}
\begin {aligned}
\rho (R{D_\alpha }(G')) - \rho (R{D_\alpha }(G))
  \geq& {\textbf{x}^{\rm T} }R{D_\alpha }(G')\textbf{x} - {\textbf{x}^{\rm T} }R{D_\alpha }(G)\textbf{x}\\
  =&\alpha \left[\frac{{{n_1} -{n_2} - 1}}{2}{x_1}^2 + \frac{{({n_1} - 1)}}{2}{x_1}^2 - \frac{{{n_2}}}{2}{x_2}^2\right]\\
  &+ (1 - \alpha )[({n_1} - 1){x_1}^2-{n_2}{x_1}{x_2} ]\\
 = &({n_1} - 1){x_1}^2- {n_2}{x_1}{x_2} - \frac{{\alpha {n_2}}}{2}{({x_1} - {x_2})^2}.
  \end {aligned}
\end{equation}
For the sake of convenience, let $\rho=\rho (R{D_\alpha }(G))$. From $R{D_\alpha }(G)\textbf{x}= \rho \textbf{x}$, we easily find that, for $1\leq i\leq \chi$,
\[\rho {x_i} = \alpha (n - \frac{{{n_i} + 1}}{2}){x_i} + (1 - \alpha )(\sum\limits_{k = 1}^\chi  {{n_k}} {x_k} - \frac{{{n_i} + 1}}{2}{x_i}).
\]
It means that
\begin{equation}\label{23}
(\rho  - \alpha n + \frac{{{n_i} + 1}}{2}){x_i} = (1 - \alpha )\sum\limits_{k = 1}^\chi  {{n_k}} {x_k}.
\end{equation}
Therefore, it follows from (\ref{23}), along with $n_{1}\geq n_{k}$, that $x_{1}\leq x_{k}$ for $k=2,3,\ldots,\chi$ and
\begin{equation}\label{24}
(\rho  - \alpha n + \frac{{{n_1} + 1}}{2}){x_1} = (\rho  - \alpha n + \frac{{{n_2} + 1}}{2}){x_2}.
\end{equation}
Now, plugging (\ref{24}) into (\ref{22}), one obtains that
\begin{equation}\label{25}
\begin {aligned}
\rho (R{D_\alpha }(G')) - \rho (R{D_\alpha }(G))
&\geq({n_1} - 1){x_1}^2- {n_2}{x_1}{x_2} - \frac{{\alpha {n_2}}}{2}{({x_1} - {x_2})^2}\\
  &={x_1}^2\left[ {({n_1} - 1) - {n_2}\frac{{{x_2}}}{{{x_1}}} - \frac{{\alpha {n_2}}}{2}{{(1 - \frac{{{x_2}}}{{{x_1}}})}^2}} \right]\\
  &= \frac{{{x_1}^2}}{{\rho  - \alpha n + \frac{{{n_2} + 1}}{2}}}\left[({n_1} - 1 - {n_2})(\rho  - \alpha n + \frac{{{n_1} + 1}}{2}) \right.\\
   &\;\;\;\;\left. - \frac{{({n_1} - 1)({n_1} - {n_2})}}{2} - \frac{{\alpha {n_2}}}{8}\frac{{{{({n_2} - {n_1})}^2}}}{{\rho  - \alpha n + \frac{{{n_2} + 1}}{2}}}\right].
  \end {aligned}
\end{equation}
As we all know,
\begin{equation}\label{26}
\rho (R{D_\alpha }(G))\geq\mathop {\min }\limits_{1 \leq i \leq n} (\sum\limits_{j = 1}^n {R{D_\alpha }{{(G)}_{i,j}}} ) = n - \frac{{{n_1} + 1}}{2},
\end{equation}
where ${R{D_\alpha }{{(G)}_{i,j}}}$ is the $(i,j)$th entry of $R{D_\alpha }(G)$. According to (\ref{26}) and $n_1>n_2+1$, we have $\rho  - \alpha n + \frac{{{n_2} + 1}}{2} > 2(1 - \alpha ){n_2}>2\alpha {n_2}$ for $0\leq\alpha \leq \frac{7}{16}$, which implies that
\begin{equation*}
\begin {aligned}
  ({n_1} - 1 - {n_2})&(\rho  - \alpha n + \frac{{{n_1} + 1}}{2}) - \frac{{({n_1} - 1)({n_1} - {n_2})}}{2} - \frac{{\alpha {n_2}}}{8}\frac{{{{({n_2} - {n_1})}^2}}}{{\rho  - \alpha n + \frac{{{n_2} + 1}}{2}}} \hfill \\
  &> ({n_1} - 1 - {n_2})(\rho  - \alpha n + \frac{{{n_1} + 1}}{2}) - \frac{{({n_1} - 1)({n_1} - {n_2})}}{2} - \frac{{{{({n_2} - {n_1})}^2}}}{16} \hfill \\
  & = ({n_1} - 1 - {n_2})(\rho  - \alpha n + \frac{1}{2} - \frac{{{n_1} + 1 - {n_2}}}{16}) - \frac{{{n_2}}}{2} - \frac{{1  }}{16} \hfill \\
  &=({n_1} - 1 - {n_2})\left[(1 - \alpha )\sum\limits_{k = 1}^\chi  {{n_k}} \frac{{{x_k}}}{{{x_1}}} - \frac{{{n_1}}}{2} - \frac{{{n_1} + 1 - {n_2}}}{16}\right] - \frac{{{n_2}}}{2} - \frac{1}{16}\;\;\;\;\;\;(\text{by}\;\; (\ref{23})) \hfill \\
  & \geq ({n_1} - 1 - {n_2})\left[(1 - \alpha )\sum\limits_{k = 3}^\chi  {{n_k}}  + (\frac{7}{16} - \alpha )(n_1 +n_2) +\frac{1}{8}(n_2-1)\right], \hfill \\
\end {aligned}
\end{equation*}
where the last inequality follows because ${n_1} - 1 - {n_2}\geq 1$ and $x_{1}\leq x_{k}$ for $k=2,3,\ldots,\chi$. Finally, it follows from this, along with (\ref{25}), that $\rho (R{D_\alpha }(G')) - \rho (R{D_\alpha }(G))>0$ for $0\leq\alpha\leq\frac{7}{16}$. This  contradicts our previous assumption. Therefore, $G=T_{n,\chi}$.
\end{proof}

A vertex subset $S\in V(G)$ is called an independent set of a graph $G$ if any two vertices of $S$ is not adjacent in $G$. The independent number of $G$ is the number of elements in the  maximum independent set of graph $G$.

\paragraph{Theorem 4.5.} Let $G$ be a simple connected undirected graph of $n$ vertices with independence number $k$ and $0\leq \alpha<1$. Then
\[
\rho(RD_{\alpha}(G))\leq \frac{(1+\alpha)n-\frac{k}{2}-\frac{3}{2}+\sqrt{((1-\alpha)n+2\alpha k-\frac{3k}{2}-\frac{1}{2})^{2}+4(1-\alpha)^{2}k(n-k)}}{2}.
\]
Furthermore, the equality holds if and only if $G=\overline{K_{k}}\vee K_{n-k}$.
\begin{proof}
Since $G$ is a connected graph with independence number $k$. Then,  Proposition 3.2 implies that  $\rho(RD_{\alpha}(G))\leq \rho(RD_{\alpha}(\overline{K_{k}}\vee K_{n-a}))$. Now, according to Corollary 2.4, the equation
\[\rho(RD_{\alpha}(\overline{K_{k}}\vee K_{n-k})=
 \frac{(1+\alpha)n-\frac{k}{2}-\frac{3}{2}+\sqrt{((1-\alpha)n+2\alpha k-\frac{3k}{2}-\frac{1}{2})^{2}+4(1-\alpha)^{2}k(n-k)}}{2}
 \]
is true.
\end{proof}

\section{Concluding remarks}

This work is mainly concerned with spectral properties of the generalized reciprocal distance matrix $RD_{\alpha}(G)$, which is the unified way of the reciprocal distance matrix $RD(G)$ and reciprocal distance signless Laplacian matrix $RQ(G)$. Several basic spectral properties of $RD_{\alpha}(G)$ and some bounds on generalized reciprocal distance spectral radius are established. We completely determine the generalized reciprocal distance spectra of some special graphs, which are used to discuss the positive semidefinite properties of $RD_{\alpha}(G)$. We also characterize the extremal graphs with maximal generalized reciprocal distance spectral radius in several kinds of simple connected graphs with precise graph invariants, vertex connectivity, edge connectivity, chromatic number, independence number and so on.

It is worth mentioning that there is a small amount of research for the new introduced matrix $RQ(G)$ (see \cite{Alhevaz,Medina}). Studying spectral properties of $RD_{\alpha}(G)$ not only generalizes the related results of $RD(G)$, but also promotes the spectral research of $RQ(G)$. For example, if we take $\alpha=0$ in Theorem 4.1, then ${K_r} \vee ({K_1} \cup {K_{n - r - 1}})$ has the maximal reciprocal distance spectral radius among connected graphs with precise vertex connectivity, which is exactly a result in \cite{Su}. If we take $\alpha=\frac{1}{2}$ in Theorem 4.1, then we arrive at:

\paragraph{Theorem 5.1.} The graph ${K_r} \vee ({K_1} \cup {K_{n - r - 1}})$ is the unique graph with maximal reciprocal distance signless Laplacian spectral radius in $G\in G_{n}^{r}$ with $1\leq r\leq n-2$.\\
\\
However, that doesn't always seem to happen. For example, it is proved in our Theorem 4.4 that the Tur\'{a}n graph $T_{n,\chi}$ has maximal spectral radius of $RD_{\alpha}(G)$ in $G\in \widetilde{G}_{n}^{\chi}$ with $0\leq \alpha \leq \frac{7}{16}$. We don't know if the result holds for $\alpha=\frac{1}{2}$.
 This observation raises the following problem:
\begin{enumerate}\setcounter{enumi}{0}
  \item Whether there exists $\alpha_0\in [\frac{1}{2}, 1)$ such that $T_{n,\chi}$ has maximal spectral radius of $RD_{\alpha}(G)$ in $G\in \widetilde{G}_{n}^{\chi}$ for $0\leq \alpha\leq \alpha_{0}$, or determine the graphs with maximal spectral radius of $RD_{\alpha}(G)$ in $G\in \widetilde{G}_{n}^{\chi}$ for $\alpha\in (\frac{7}{16}, 1)$.
\end{enumerate}

There are, of course, many other problems to be considered in the future work, for instance:
\begin{enumerate}\setcounter{enumi}{1}
  \item Characterize the graphs with  minimal or maximal spectral radius of $RD_{\alpha}(G)$ in several kinds  of graphs, such as trees, unicyclic graphs, the graphs with given matching number and so on.
  \item Given some special graph classes, determine the smallest $\alpha_{0}\in(0,\frac{1}{2}]$ such that $R{D_\alpha }(G)$ is positive semidefinite whenever $\alpha_{0}\leq \alpha\leq 1$.
\end{enumerate}

Finally, we point out that, the generalized reciprocal distance energy  for a connected graph $G$ of order $n$ can be defined as
\[
E_{RD_\alpha}(G)=\sum\limits_{i = 1}^n\left|\lambda_{i}(RD_{\alpha}(G))-\frac{2\alpha H(G)}{n}\right|.
\]
It is clear that $E_{RD_0}(G)$ is the reciprocal distance energy of $G$ \cite{Gungor}, while $E_{RD_{1/2}}(G)$ is exactly the reciprocal distance
signless Laplacian energy of $G$ \cite{Alhevaz,Medina2}, Therefore, it is entirely possible for us to focus on research these energies in a unified way.

\end{document}